\documentclass{article}
\usepackage{amssymb}
\usepackage{latexsym}
\usepackage{graphicx}

\def \ee { {2^{-n}}}
\def \eee{\cdot 2^{-n}}
\def \ddd {,\ldots, }
\def \lp {\mathcal L(\overline p)}

\def \IQ {\mathbb Q}
\def \IR {\mathbb R}
\def \IN {\mathbb N}

\def \proof{\noindent {\bf Proof\ \  }}

\def \dom{{\rm dom}}
\def \range {{\rm range}}
\def \s  {\Sigma^*}
\def \om {\Sigma^\omega}
\def \mod {{\rm mod}}
\def \an {\mbox{\ \ and\ \ }}

\def \In{\subseteq}
\def \pf{:\hspace{0.6ex}\subseteq \hspace{-0.4ex}}

\newcommand{\qq}{\phantom{.}\hfill$ \Box$}

\newcommand{\mto}{\rightrightarrows}

\newcommand{\mmto}{\mbox{
\setlength{\unitlength}{1em}
\begin{picture}(0.4,0)
\makebox(0,0.6){$\mbox{\scriptsize \raisebox{0.083em}{$|$}}   \hspace*{-1.1ex}\mto$}
\end{picture}
}}

\newtheorem{definition}{Definition}[section]{\bf}{\it}
\newtheorem{theorem}[definition]{Theorem}{\bf}{\it}
{\bf}{\it}
\newtheorem{lemma}[definition]{Lemma}{\bf}{\it}
\newtheorem{corollary}[definition]{Corollary}{\bf}{\it}
\newtheorem{proposition}[definition]{Proposition}{\bf}{\it}
{\bf}{\it}
{\bf}{\it}
{\bf}{\it}

\begin{document}

\author{Klaus Weihrauch\\ \\ FernUniversit\"at in Hagen}
\title{Intersection points of planar curves can be computed}
\date{}

\maketitle

\begin{abstract} Consider two paths $\phi,\psi:[0;1]\to [0;1]^2$ in the unit square such that $\phi(0)=(0,0)$, $\phi(1)=(1,1)$, $\psi(0)=(0,1)$ and $\psi(1)=(1,0)$. By continuity of $\phi$ and $\psi$ there is a point of intersection. We prove that from $\phi$ and $\psi$ we can compute closed intervals $S_\phi,S_\psi \In [0;1]$ such that $\phi(S_\phi)=\psi(S_\psi)$.

\end{abstract}

\section{Introduction}\label{seca}A {\em path} in the Euclidean plane is a continuous function $f: [0;1]\to \IR^2$, a {\em curve} is the range of a path. The following is known about planar curves.

\begin{theorem}\label{t1}
Let  $\phi,\psi :[0;1]\to [0;1]^2$  be two paths in the unit square such that
\begin{eqnarray}
\label{f1}& \phi(0)=(0,0), \ \phi(1)=(1,1), \  \psi(0)=(0,1)\an   \psi(1)=(1,0) \,.
\end{eqnarray}
Then the two curves $\range(\phi)$ and $\range(\psi)$ intersect.
\end{theorem}

\begin{figure}[htbp]
\setlength{\unitlength}{.2ex}
%\linethickness{0.7pt}
\begin{picture}(200,165)(-210,-25)

\put(-130,0){\vector (1,0){360}}
\put(0,100){\line (1,0){100}}
\put(0,-30){\vector (0,1){160}}
\put(100,0){\line (0,1){100}}
\put(0,0){\circle*{6}}
\put(243,0){\makebox(0,0)[cc]{$x_1$}}
\put(0,138){\makebox(0,0)[cc]{$x_2$}}

\thicklines
\qbezier(0,0)(20,80)(100,100)
\qbezier(0,100)(50,10)(100,0)

\put(-100,0){\line(1,0){100}}
\put(100,0){\line(1,0){100}}
\put(-100,100){\line(1,0){100}}
\put(100,100){\line(1,0){100}}

\put(54,70){\makebox(0,0)[cc]{$\phi$}}
\put(40,30){\makebox(0,0)[cc]{$\psi$}}

\put(-50,90){\makebox(0,0)[cc]{$g$}}
\put(-50,10){\makebox(0,0)[cc]{$f$}}
\put(150,90){\makebox(0,0)[cc]{$f$}}
\put(150,10){\makebox(0,0)[cc]{$g$}}

\put(13,-10){\makebox(0,0)[cc]{$\phi(0)$}}
\put(13,108){\makebox(0,0)[cc]{$\psi(0)$}}
\put(100,-10){\makebox(0,0)[cc]{$\psi(1)$}}
\put(100,110){\makebox(0,0)[cc]{$\phi(1)$}}

\put(-100,-10){\makebox(0,0)[cc]{$f(-1)$}}
\put(-100,110){\makebox(0,0)[cc]{$g(-1)$}}
\put(200,-10){\makebox(0,0)[cc]{$g(2)$}}
\put(200,110){\makebox(0,0)[cc]{$f(2)$}}

\end{picture}
\caption{Intersecting curves $\phi$ and $\psi$ with extensions $f$ and $g$.} \label{fig1}
\end{figure}
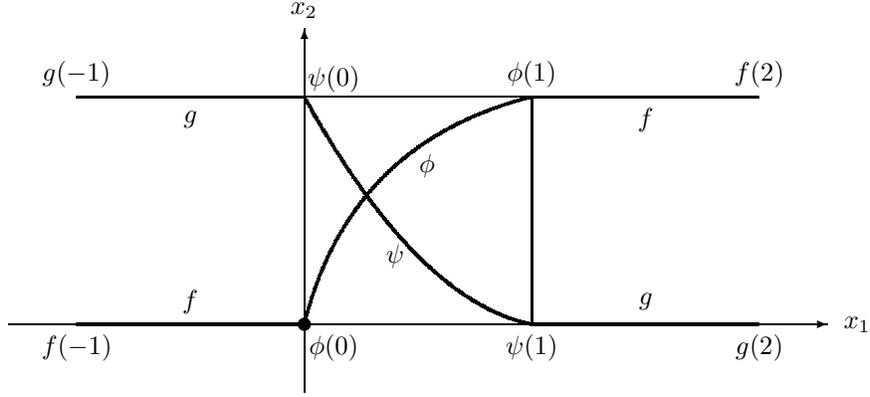

Figure~\ref{fig1} visualizes the theorem. In Markov-style computable analysis \cite{Kus84}
Manukyan \cite{Man76a} has proved a related theorem (in Russian), cited in \cite[Page 279]{Kus99} as follows:

\begin{theorem}[Manukyan]\label{t2a}There are two constructive (and therefore continuous) planar curves  $ \varphi_1$ and $ \varphi_2$ such that
\begin{eqnarray}
\label{f5aa}& \varphi_1(0)=(0,0), \  \varphi_1(1)=(1,1),\
 \varphi_2(0)=(0,1), \   \varphi_1(1)=(1,0)\,, \\
\label{f6aa}&\mbox{for every $0<t< 1$ both $ \varphi_1(t)$  and $ \varphi_2(t)$ belong to the open unit square,}\\
\label{f7aa}&\mbox{the paths of $ \varphi_1$ and $ \varphi_2$ do not intersect.}
\end{eqnarray}
\end{theorem}
While in (Grzegorczyk-Lacombe-  \cite{Grz55,Lac55a}) computable analysis the following has been proved \cite{Wei19}:

\begin{theorem}[Weihrauch]\label{t1a}
If $\phi$ and $\psi$ in Theorem~\ref{t1} are computable then there is a computable point $x\in\range(\phi)\cap\range(\psi)$.
\end{theorem}

This is not contradictory. In Markov's approach only functions on the computable real numbers which are encoded by G{\"o}del numbers are considered and computations transform G{\"o}del numbers to G{\"o}del numbers. While in the Grzegorczyk-Lacombe-approach all real numbers are considered, where real numbers are encoded by (fast converging) Cauchy sequences of rational numbers and computations transform infinite Cauchy sequences to infinite Cauchy sequences.

In this article we prove Claim~6.2 from \cite{Wei19}:
\begin{theorem} \label{t2}
Let $\mathcal T$ be the multi-valued operator  mapping every pair
$\phi,\psi:[0;1]\to [0;1]^2$  of paths in the unit square such that
\begin{eqnarray}
\label{f2}& \phi(0)=(0,0), \ \phi(1)=(1,1), \  \psi(0)=(0,1)\an   \psi(1)=(1,0)
\end{eqnarray}
to some pair $(S_\phi,S_\psi)$ of  closed intervals such that $\phi(S_\phi)=\psi(S_\psi)$. Then the operator  $\mathcal T$ is computable.
\end{theorem}

 Theorem~\ref{t1a} follows straightforwardly from Theorem~\ref{t2}.
In the proof from $\phi$ and $\psi$ we compute sequences  $I_0\supseteq I_1\supseteq I_2 \supseteq\ldots$ and  $J_0\supseteq J_1\supseteq J_2 \supseteq\ldots$ of closed intervals with rational endpoints such that
$\phi(\bigcap_iI_i)=\psi(\bigcap _iJ_i)$.

Curves (even computable ones) can be much more complicated than the examples shown in Figure~\ref{fig1}.  Consider, for example, space-filling curves or curves with infinitely many spirals, each of which containing infinitely many sub-spirals etc. infinitely often or curves with ``completely" chaotic behavior.

This article is a contribution to computable analysis. There are various non equivalent definitions of computability in analysis. One of these is Markov's constructive analysis \cite{Kus84,Kus99}. Theorem~\ref{t2a} is a result in this theory.
We use "TTE", an approach which is based on ideas from \cite{Grz55,Grz57,Lac55a}.
In TTE computability
on $\{0,1\}^*$ and  Cantor space $\{ 0,1\}^\omega$ (the finite and infinite $0$-$1$-sequences)
is defined explicitly (e.g. by Turing machines  with finite or infinite one-way input and output tapes) and computability on other sets $X$ is induced via representations $\delta \pf \{0,1\}^*\to X$ or
$\delta\pf \{ 0,1\}^\omega \to X$ (partial surjective) where finite or infinite $0$-$1$-sequences  are interpreted as names and computations are performed on names. We consider canonical representations of
the real numbers, open subsets, closed subsets, compact subsets and real functions. Equivalently any finite alphabet  $\Sigma$  (with at least two elements) can be used instead of $\{0,1\}$.
We assume that the reader is familiar with the basic concepts of TTE. Details can be found in \cite{Wei00,BHW08,WG09}.
\medskip

For technical reasons we extend $\phi$ and $\psi$ trivially to continuous functions $f,g:[-1;2]\to\IR^2$  which intersect in the same way as $\phi$ and $\psi$, that is, $\phi(s)=\psi(t)\iff f(s)=g(t)$ (see Figure~\ref{f1}):
\begin{eqnarray}
\label{f15}
f(t)& := &
\left\{\begin{array}{rcl}
(t,0)& \mbox{if}& -1\leq t\leq 0\\
\phi(t) & \mbox{if}& 0\leq t\leq 1\\
(t,1) & \mbox{if}& 1\leq t\leq2 \,,
\end{array}\right.\\
\label{f14}
g(t)& := &
\left\{\begin{array}{rcl}
(t.1)& \mbox{if}& -1\leq \leq 0\\
\psi(t) & \mbox{if}& 0\leq t\leq 1\\
(t,0) & \mbox{if}& 1\leq t\leq 2\,,
\end{array}\right.
\end{eqnarray}

We will consider closed  rational sub-intervals $I=[a_I;b_I]$ and $J=[a_J;b_J]$  of the real interval $[-1;2]$
such that for the restrictions $f|_I$ of $f$ to $I$ and $g|_J$ of $g$ to $J$,
\begin{eqnarray*}
&&\mbox{the end-points $f|_I(a_I)$ and $f|_I(b_I)$ are not in $g(J)$}\an\\
&&\mbox{the end-points $g|_J(a_J)$ and $g|_J(b_J)$ are not in $f(I)$}\,.
\end{eqnarray*}
Since $f(I)$ and $g(J)$ are compact this means
\begin{eqnarray}\label{f8}
\alpha_{IJ}:= \min(d_s(\{f(a_I),f(b_I)\},g(J)), d_s(\{g(a_J),g(b_J)\},f(I)))>0\,.
\end{eqnarray}
where
$d_s(A_1,A_2)  :=  \inf\{\|z_1-z_2\|\mid z_1\in A_1,z_2\in A_2\}$.

We will approximate  $f|_I$  and $g|_J$   by rational polygon paths
$h$ and $h'$, respectively, and consider the intersections  $(s,t)$, that is, pairs such that $h(s)=h'(t)$. In order to keep this number finite we consider only pairs $(h,h')$ such that $\range(h)\cap \range(h')$ contains no straight  line segment.
For such pairs every intersection $(s,t)$ is either a crossing or tangent. As a central lemma we will prove that the parity (even or odd) of the  number of crossings does not depend on $h$ and $h'$ (it is an invariant of $(f|_I,g|_J)$). We call it the crossing parity of the pair $(f|_I ,g|_J)$.
In the proof we will apply transformations of polygon paths which may change the number of crossings but do not change the parity (even or odd) of the number of crossings.

\section{The crossing parity}\label{secb}

For points $x,y\in\IR^2$ such that $x\neq y$ let $\overline{xy}\In \IR^2$ be the straight line segment from $x$ to $y$. In this article $\overline{xx}$ is not a straight line segment.

\begin{definition}\label{d1}\item
\begin{enumerate}
\item
A track is a sequence $p=((s_0,x_0),(s_1,x_1),\ldots, (s_k,x_k))$ such that  $s_i<s_{i+1}$  and $x_i\neq x_{i+1}$ for $0\leq i<k$. The points $x_0,x_1,\ldots,x_k$ are the vertices of $p$.
\item
The track $p$ spans a (polygon) path $h_p:[s_0;s_k]\to \IR^2$ by
   \begin{eqnarray}
\label{f5} h_p(s)&=&x_i+\frac{s-s_i}{s_{i+1}-s_i}(x_{i+1}-x_i)\,
 \mbox{ if } \,s_i\leq s\leq s_{i+1}\,.
\end{eqnarray}
\end{enumerate}
\end{definition}
\noindent  By (\ref{f5}), $h_p(s_i)=x_i$ and $h_p[s_i;s_{i+1}]=\overline{x_ix_{i+1}}$.
\noindent

For tracks $p,q$ we want to count the number of crossings of the paths $h_p$ and $h_q$. In order to keep this number finite we consider only pairs $(p,q)$  such that $\range(h_p)$ and $\range(h_q)$ have no common straight line segment.
 Furthermore, we will not count all intersections of $h_p$ and $h_q$ but only  ``proper" crossings.

\begin{definition}\label{d3}
Let $(p,q)$ where $p=((s_0,x_0),s_1,x_1),\ldots, (s_k,x_k))$ and  $q=((t_0,y_0),t_1,y_1),\ldots, (t_l,y_l))$
be a pair of tracks such that $\range(h_p)$ and $\range(h_q)$ have no common  straight line segment.
\begin{enumerate}
\item
An intersection of $p$ and $q$ is a pair $(s,t)$ such that
$s_0< s<s_k$,  $t_0< t<t_l$ and
$h_p(s) = h_q(t)$.
We call  x:= $h_p(s) = h_q(t)$ the corresponding intersection point.

\item For an intersection $(s,t)$  of $p$ and $q$  with intersection point $h_p(s)=h_q(t)=x$ let $\delta_{st}>0 $ be a number  such that $B(x,\delta_{st})\setminus\{x\}$ contains no vertex of~$p$ and no vertex of~$q$. Let
\begin{eqnarray*}
&s_< := \inf\{s'<s\mid h_p[s';s]\In B(x,\delta_{st})\},& x_<:=h_p(s_<)\,,\\
&s_> := \sup\{s'>s\mid h_p[s;s']\In B(x,\delta_{st})\},& x_>:=h_p(s_>)\,,\\
&t_< := \inf\{t'<t\mid h_q[t';t]\In B(x,\delta_{st})\},& y_<:=h_q(t_<)\,,\\
&t_> := \sup\{t'>t\mid h_q[t;t']\In B(x,\delta_{st})\},& y_>:=h_q(t_>)\,.
\end{eqnarray*}
If on the boundary of $B(x,\delta_{st})$ the four points $x_<,\;x_>,\;y_<$ and $y_>$  occur  in the order
$(x_<,y_<,x_>,y_>)$ or in the order $(x_<,y_>,x_>,y_<)$
\footnote{that is, on the boundary of $B(x,\delta_{st})$  the four points alternate in $x$ and $y$}, we call $(s,t)$ a crossing and $x$ the corresponding crossing point, else $x$ is a touch point.

\item Let ${\rm CN}(p,q)$ the number of crossings of $p$ and $q$ and let \\ $\pi(p,q):= {\rm CN}(p,q)\; \mod\; 2$ be  its parity ($\:0=$ even and $1=$ odd).
\end{enumerate}
\end{definition}
Obviously,  \ $s_< <s_i<s$  for no number $i$, \ $x_<=h_p(s_<)\in \partial B(x,\delta_{st})$ and $h_p[s_<;s]=\overline {x_< x}$ where $\partial A$ denotes the boundary of $A\In\IR^2$. This is true correspondingly for $s_>,t_<$ and $t_>$.

\begin{figure}[htbp]
\setlength{\unitlength}{2.0pt}
%\linethickness{0.7pt}
\begin{picture}(400,42)(-34,-4)

\newsavebox{\ballu}
\savebox{\ballu}{
\qbezier(15,0)(15,6.21)(10.61,10.61)
\qbezier(10.61,10.61)(6.21,15)(0,15)
\qbezier(15,0)(15,-6.21)(10.61,-10.61)
\qbezier(10.61,-10.61)(6.21,-15)(0,-15)
\qbezier(-15,0)(-15,6.21)(-10.61,10.61)
\qbezier(-10.61,10.61)(-6.21,15)(0,15)
\qbezier(-15,0)(-15,-6.21)(-10.61,-10.61)
\qbezier(-10.61,-10.61)(-6.21,-15)(0,-15)
\put(0,0){\circle*{1.5}}
}

\put(20,20){\usebox{\ballu}}
\put(60,20){\usebox{\ballu}}
\put(100,20){\usebox{\ballu}}
\put(140,20){\usebox{\ballu}}

\thinlines
\qbezier(9.39,9.39)(9.39,9.39)(30.61,30.61)
\thicklines
\qbezier(5,20)(5,20)(35,20)
\thinlines
\put(1,20){\makebox(0,0)[cc]{$x_>$}}
\put(40,20){\makebox(0,0)[cc]{$x_<$}}
\put(8,6){\makebox(0,0)[cc]{$y_<$}}
\put(33,34){\makebox(0,0)[cc]{$y_>$}}

\put(60,20){\line(2,1){13.55}}
\put(60,20){\line(-1,-1){10.5}}
\thicklines
\qbezier(60,20)(60,20)(49,30)
\qbezier(60,20)(60,20)(74,15)
\thinlines
\put(47,32){\makebox(0,0)[cc]{$x_<$}}
\put(79,14){\makebox(0,0)[cc]{$x_>$}}
\put(47,8){\makebox(0,0)[cc]{$y_<$}}
\put(78,27){\makebox(0,0)[cc]{$y_>$}}

\put(100,20){\line(2,1){13.55}}
\qbezier(100,20)(100,20)(89,30)
\thicklines
\qbezier(100,20)(100,20)(95,6)
\qbezier(100,20)(100,20)(114,15)
\thinlines

\put(86,32){\makebox(0,0)[cc]{$y_>$}}
\put(117,29){\makebox(0,0)[cc]{$y_<$}}
\put(94,3,5){\makebox(0,0)[cc]{$x_>$}}
\put(119,14){\makebox(0,0)[cc]{$x_<$}}

\put(140,19.6){\line(1,0){15}}
\put(140,20.4){\line(1,0){15}}
\thicklines
\qbezier(140,20)(140,20)(129,30)
\qbezier(140,20)(140,20)(151,30)
\thinlines
\put(158,17){\makebox(0,0)[cc]{$y_<=y_>$}}
\put(127,32){\makebox(0,0)[cc]{$x_<$}}
\put(154,32){\makebox(0,0)[cc]{$x_>$}}

\put(20,-4){\makebox(0,0)[cc]{(a)}}
\put(60,-4){\makebox(0,0)[cc]{(b)}}
\put(100,-4){\makebox(0,0)[cc]{(c)}}
\put(140,-4){\makebox(0,0)[cc]{(d)}}

\end{picture}
\caption{$h_p$ and $h_q$ in the ball $B(x,\delta_{st})$ for $x=h_p(s)=h_q(t)$.} \label{fig4}
\end{figure}
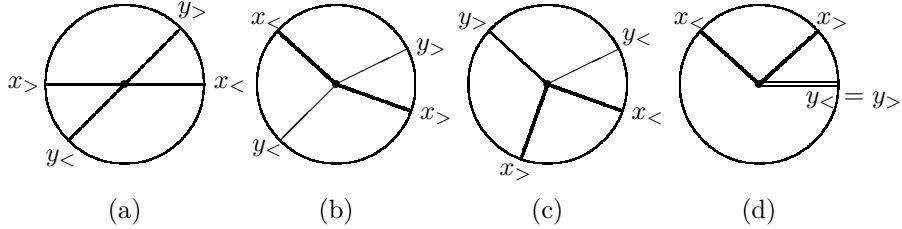

Figure~\ref{fig4} shows several kinds of intersection of $p$ and $q$
(thin lines for $h_q$ and thick lines for $h_p$ ) the first two of which are crossings. In (a) possibly the center $x$ is no vertex of $p$ or no vertex  of $q$.

The definition of a crossing $(s,t)$ depends on the number $\delta_{st}$ only formally.

\begin{lemma}\label{l2}
The definition of a crossing does not depend on the choice of $\delta_{st}$.
\end{lemma}

\proof
Let $(s,t)$ be a crossing of $p$ and $q$ defined via some $\delta_{st}$.
Let $0<\overline \delta <\delta_{st}$.

Let $\overline s_< := \inf\{s'<s\mid h_p[s';s]\In B(x,\overline \delta)\}$ and $\overline x_<:=h_p(\overline s_<)$. Then
$\overline x_< =\overline {x_<x}\cap \delta B(x,\overline \delta)$.
This is true correspondingly for the other three cases. Obviously the four points on $B(x,\overline \delta)$ alternate in the same way as the four corresponding points on $B(x, \delta_{st})$.
\qq\\

Notice that $\{x_<,x_>\}\cap \{y_<,y_>\}=\emptyset$ since $\range(h_p)$ and $\range(h_q)$ have no common  straight line segment.
In the applications below the endpoints  of $h_p$ are not in $\range(h_q)$ and  the endpoints  of $h_q$ are not in $\range(h_p)$. Therefore it suffices to consider  $s_0< s<s_k$ and $t_0< t<t_l$ in the definition of intersections.

We introduce a separation concept for tracks  $p$ and $q$ which induces that $\range(h_p)$ and $\range(h_q)$ have no common  straight line segment.

\begin{definition}\label{d6}
Let $p=((s_0,x_0),(s_1,x_1),\ldots, (s_k,x_k))$ be a track. \begin{enumerate}
 \item Define  $\mathcal V(p):=\{x_0,x_1,\ldots,x_k\}$\,.
 \item For $x,y\in\IR^2$ with $x\neq y$ let $l(x,y)$ be the straight line through  $x$ and~$y$.

\item  Define
\begin{eqnarray}
 \label{f1a}\mathcal {L}(p) &:=&\bigcup \{l(x_{i-1},x_i)\mid
 1\leq i\leq k\}\,,\\
\label{f1b} \overline \mathcal L (p) & :=& \bigcup\{l(x_i,x_j)\mid 0\leq i<j\leq k,\ x_i\neq x_j\}\,.
\end{eqnarray}
 \item We call tracks $p$ and $q$ weakly separated,
 $p\bowtie q$, iff
\begin{eqnarray}
 \label{f12}
\mathcal V(p)\cap \mathcal{L}(q)&=&\emptyset \an\\
  \label{f13}
\mathcal V(q)\cap \mathcal{L}(p)&=&\emptyset \,.
 \end{eqnarray}
\end{enumerate}
\end{definition}

\begin{figure}[htbp]
\setlength{\unitlength}{1pt}

\begin{picture}(200,0)(-80,100)
\put(0,30){\line(2,3){46.5}}
\put(20,0){\line(1,1){80}}
\put(20,100){\line(2,-3){67}}
%
%\qbezier(0,99)(50.4,0)(50.4,0)
%\qbezier(0,70)(0,70)(100,20)
%\qbezier(42,0)(42,0)(31,100)

\thicklines
\put(20.23,60){\circle*{4}}
\put(33.2,80){\circle*{4}}
\put(59.5,40){\circle*{4}}
\put(40,19.6){\circle*{4}}

\qbezier(20.23,60)(20.23,60)(33.2,80)
\qbezier(33.2,80)(59.5,40)(59.5,40)
\qbezier(59.5,40)(40,19.6)(40,19.6)

\thinlines
\put(-20,80){\line(4,-1){140}}
\put(34.5,66.2){\circle{4}}
\put(34,60){\makebox(0,0)[cc]{$x$}}

\end{picture}

\begin{picture}(200,95)(-280,4)

\put(0,30){\line(2,3){46.5}}
\put(20,0){\line(1,1){80}}
\put(20,100){\line(2,-3){67}}

\qbezier(0,99)(50.4,0)(50.4,0)
\qbezier(0,70)(0,70)(100,20)
\qbezier(42,0)(42,0)(31,100)

\linethickness{0.3ex}

\put(20.23,60){\circle*{4}}
\put(33.2,80){\circle*{4}}
\put(59.5,40){\circle*{4}}
\put(40,19.6){\circle*{4}}

\qbezier(20.23,60)(20.23,60)(33.2,80)
\qbezier(33.2,80)(59.5,40)(59.5,40)
\qbezier(59.5,40)(40,19.6)(40,19.6)

\thinlines
\put(-20,80){\line(4,-1){140}}
\put(65,58.86){\circle{4}}
\put(65,66){\makebox(0,0)[cc]{$y$}}
\end{picture}

\caption{$\mathcal L(p)$ and $\mathcal{\overline L}(q)$}
  \label{fig2}
\end{figure}
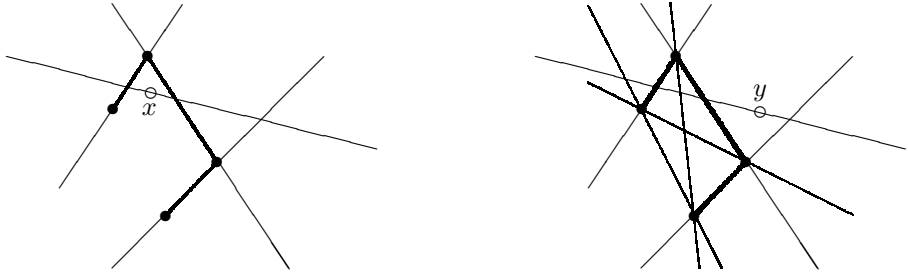

Figure~\ref{fig2} shows on the left the set $\mathcal{L}(p)$ of a track $p$ and a straight line through a point $x\not\in \mathcal L(p)$ and on the right the set $\overline \mathcal L(q)$ of a track $q$ and a straight line through a point $y\not\in \overline\mathcal  L(q)$.

\begin{lemma}\label{l7}$ $
\begin{enumerate}
 \item
 If $y\not \in \mathcal{L}(p)$ then every straight line through the point $y$  intersects every straight line from $\mathcal L(p)$ at most once.
\item If $\,\mathcal V(p)\cap \mathcal{L}(q)=\emptyset$ or
$\mathcal V(q)\cap \mathcal{L}(p)=\emptyset$ then then $\range(h_p)$ and $\range(h_q)$ have no common  straight line segment.
\end{enumerate}
\end{lemma}

\proof
If a straight line through $y$  intersects a straight line from $ \mathcal{L}(p)$ twice then  $y\in  \mathcal{L}(p)$. Contradiction.

Let $q=((t_0,y_0) \ddd (t_m,y_m))$. Suppose for some $i$ and $j$,  $\overline {x_ix_{i+1}}$ and
$\overline {y_jy_{j+1}}\ $ have a  common  straight line segment. Then $x_i\in l(y_j,y_{j+1})\In \mathcal L(q)$, but
$\mathcal V(p)\cap \mathcal{L}(q) =\emptyset$. Correspondingly,
$\mathcal V(q)\cap \mathcal{L}(p)\not=\emptyset$.
\qq\\

As an essential tool we will use local transformations of tracks which leave the crossing parity invariant. The following lemma justifies these transformations.

\begin{lemma}\label{l8}
Let
\begin{eqnarray*}
 p& = &((s_0,x_0),(s_1,x_1),\ldots, (s_k,x_k))\,, \\
 q_1& =&((r,y),(t_1,z_1),(r',y')) \an \\
 q_2& =&((r,y),(t_2,z_2),(r',y'))
\end{eqnarray*}
  be tracks and let $B$ be a ball such that
\begin{eqnarray}
\label{f30}&\{y,y'\}\cap \mathcal L(p)=\emptyset\,,\\
\label{f31}&\{y,y',z_1,z_2\}\In B\,,\\
\label{f32}&\{x_0,x_k\}\cap B=\emptyset\,.
\end{eqnarray}
Then   $\pi(p,q_1)= \pi(p,q_2)$.
\end{lemma}

\proof
Remember that by Definition~\ref{d1}, $r<t_1<r'$, $r<t_2<r'$,
$z_1\not\in \{y,y'\}$ and $z_2\not\in \{y,y'\}$.
By (\ref{f30}), $y,y'\not\in \range(h_p)$.

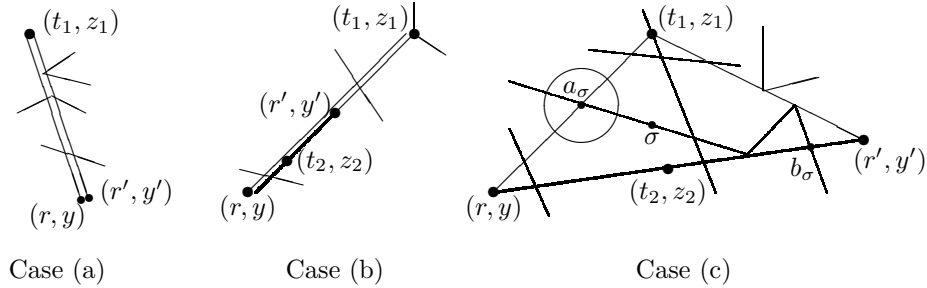
\begin{figure}[htbp]
\setlength{\unitlength}{3pt}
%\linethickness{0.7pt}
\begin{picture}(200,34)(-40,9)

\put(-18,40){\line(1,-3){7}}
\put(-19,40){\line(1,-3){7}}

\put(-18.5,40){\circle*{1.3}}
\put(-12,19){\circle*{1.}}
\put(-11,19.3){\circle*{1.}}

\put(-12,41.5){\makebox(0,0)[cc]{$(t_1,z_1)$}}
\put(-15,17){\makebox(0,0)[cc]{$(r,y)$}}
\put(-5,20){\makebox(0,0)[cc]{$(r',y')$}}

\put(-20,30){\line(2,1){4.15}}
\put(-15.7,32.15){\line(2,-1){4.3}}

\put(-16.8,35){\line(4,-1){6}}
\put(-16.8,35){\line(3,2){4}}

\put(-17,26){\line(3,-1){8}}
%\put(-22,41.7){\line(2,-1){10}}

%%%%%%%%%%%%%%%%%%%%%%%%%%%%%%%%%%%%%%%%%%%%%%%%%%%%%%%%%%5

\put(9,20){\line(1,1){20}}
\put(20,30){\line(1,1){10}}

\put(9,20){\circle*{1.3}}
\put(20,30){\circle*{1.3}}
\put(30,40){\circle*{1.3}}
\thicklines

\qbezier(10,20)(10,20)(20,30)
\put(14,24){\circle*{1.3}}

\thinlines
\put(24.5,42){\makebox(0,0)[cc]{$(t_1,z_1)$}}
\put(20,24){\makebox(0,0)[cc]{$(t_2,z_2)$}}

\put(9,18){\makebox(0,0)[cc]{$(r,y)$}}
\put(15,31){\makebox(0,0)[cc]{$(r',y')$}}

\put(8,23){\line(4,-1){8}}
\put(20,38){\line(2,-3){6}}
\put(30,44){\line(0,-1){4}}
\put(30,40){\line(3,-2){4}}

%%%%%%%%%%%%%%%%%%%%%%%%%%%%%%%

\put(40,20){\circle*{1.3}}
\put(60,40){\circle*{1.3}}
\put(86.7,26.7){\circle*{1.3}}

\put(40,20){\line(1,1){20}}
\put(60,40){\line(2,-1){26.7}}

\put(40,18){\makebox(0,0)[cc]{$(r,y)$}}
\put(90,24){\makebox(0,0)[cc]{$(r',y')$}}
\put(62,20){\makebox(0,0)[cc]{$(t_2,z_2)$}}
\put(65.5,42){\makebox(0,0)[cc]{$(t_1,z_1)$}}

\put(80,25.7){\circle*{.9}}
\put(79,23){\makebox(0,0)[cc]{$b_\sigma$}}

\put(60,28.5){\circle*{.9}}
\put(60,26.8){\makebox(0,0)[cc]{$\sigma$}}

\put(51.1,31){\circle*{.9}}
\put(51.1,31){\circle{9}}

\put(50.5,33){\makebox(0,0)[cc]{$a_\sigma$}}

\qbezier(42,28)(47,17)(47,17)

%\qbezier(50,17)(50,17))(52,32)
%\qbezier(52,32)(52,32)(70,27)

\qbezier(72,24.8)(72,24,8)(42,34)
\qbezier(72,24.8)(78,31)(78,31)
\qbezier(58.8,43)(68,20)(68,20)

%%%%%%%%%%%%%%%%%%%%%%%%%%%%%%%%%%%%%%%%%%%%%%
\qbezier(78,31)(78,31)(82,20)
\qbezier(52,38)(52,38)(71,36)

\put(74,32.85){\line(0,1){8}}
\put(74,32.85){\line(4,1){7}}

\thicklines
\qbezier(40,20)(40,20)(86.7,26.7)
\put(62,23){\circle*{1.3}}

%\put(45,25){\circle{6}}

%%%%%%%%%%%%%%%%%%%%%%%%%%%%%%%%%%%%%%%%%%%%%%%%%%%
%%%%%%%%%%%%%%%%%%%%%%%%%%%%%%%%%%%%%%%%%%%%%%%%%%%
\put(-15,10){\makebox(0,0)[cc]{Case (a)}}
\put(20,10){\makebox(0,0)[cc]{Case (b)}}
\put(64,10){\makebox(0,0)[cc]{Case (c)}}

\end{picture}
\caption{Illustration for Lemma~\ref{l8}.} \label{fig5}
\end{figure}

If $\range(h_p)$ and $\range(h_{q_1})$ have a common  straight line segment then $y\in \mathcal L(p)$ or $y'\in \mathcal L(p)$, but
$\{y,y'\}\not \in \mathcal L(p)$. Therefore, ${\rm CN}(p,q_1)$
is well-defined. Correspondingly, ${\rm CN}(p,q_2)$ is well-defined.
In Figure~\ref{fig5} the ``open ended" line segments are parts of $\range(h_p)$. We distinguish several cases, see Figure~\ref{fig5})

{\bf \boldmath Case (a) $y=y'$:} Obviously, ${\rm CN}(p,q_1)$ and
${\rm CN}(p,q_2)$ are even,
hence $\pi(p,q_1)=0=\pi(p,q_2)$.

{\bf \boldmath Case (b) $y\neq y'$, $z_2\in\overline {yy'}$ and $z_1\in l(y,y')$:} In Figure~\ref{fig5}(b) the track $q_2$ is drawn in thick lines. If $z_1\in\overline {y,y'}$ then $h_{q_1}=h_{q_2}$, hence ${\rm CN}(p,q_1)={\rm CN}(p,q_2)$.
If $z_1\not\in\overline {y,y'}$ then ${\rm CN}(p,q_1)$ and ${\rm CN}(p,q_2)$ differ by an even number. In both cases $\pi(p,q_1)=0=\pi(p,q_2)$.

{\bf \boldmath Case (c) $y\neq y'$, $z_2\in\overline {yy'}$ and $z_1\not\in l(y,y')$:} In Figure~\ref{fig5}(c) the track $q_2$ is drawn in thick lines.
Let $\Delta$ be the closed triangle with boundary
$\partial\Delta:= \range(h_{q_1}) \cup \range(h_{q_2})$
($=\overline{yz_1}\cup \overline{z_1y'}\cup \overline{yy'}$)
and let $\Delta^\circ$ be its interior.

If $\overline{yz_1}$ and $\range(h_p)$ have a common  straight line segment then $y\in\mathcal L(p)$, a contradiction by (\ref{f30}). Therefore  $\overline{yz_1}$ and $\range(h_p)$ have no common  straight line segment. This is true correspondingly for
$ \overline{z_1y'}$ and $\overline{yy'}$. Therefore
\begin{eqnarray}\label{f51}
\mbox{$\partial \Delta$  and $\range(h_p)$ have no common  straight line segment.}
\end{eqnarray}
Let $s_0<\sigma< s_k $ such that $h_p(\sigma)\in \Delta^\circ$. Let $[a_\sigma;b_\sigma]$ be the longest interval such that $s_0\leq a_\sigma<\sigma <b_\sigma\leq s_k$ and $h_p[a_\sigma;b_\sigma]\In \Delta$. Figure~\ref{fig5}(c) shows an example for $\sigma$ with $a_\sigma$ and $b_\sigma$ positioned at the images under $h_p$ and $h_p(a_\sigma)\in \range(h_{q_1})$.
Obviously, $h_p(a_\sigma), h_p(b_\sigma)\in\Delta$. We show that
$h_p(a_\sigma)$ and $ h_p(b_\sigma)$ are crossing points of $p$ and $q_1$ or of $p$ and $q_2$.

Suppose  $h_p(a_\sigma)\in \range(h_{q_1})$  \\
Since $\range(h_p)$ is a chain of  straight line segments there is some $a_\sigma<s'< \sigma$ such that $\lambda:= h_p([a_\sigma;s'])$ is a straight line segment with
$h_p(a_\sigma)\in\lambda\cap \partial\Delta$ and by (\ref{f51}) $\lambda \setminus \{h_p(a_\sigma)\}\In \Delta^\circ$.

Since $\range(h_p)$ is a chain of straight line segments and $a_\sigma$ is the smallest number $a$ with $h_p[a;\sigma]\in \Delta$, there is some $s''<a_\sigma$ such that $\lambda'':=h_p[s'';a_\sigma]$ is straight line segment  with $h_p(a_\sigma)\in\lambda''\cap \partial\Delta$ and
$(\lambda''\setminus \{h_p(a_\sigma)\}) \cap\Delta=\emptyset$ by (\ref{f51}).

If we draw a sufficiently small circle around $h_p(a_\sigma)$ then (with the terminology from Definition~\ref{d3}))  the intersections $x_<$ and $x_>$ of it with $\range(h_p)$ and the intersections $y_<$ and $y_>$ of it with $\range(h_{q_1})$ alternate on this circle  in $x$ and $y$. Therefore, for every $\sigma$ such that $h_p(\sigma)\in\Delta^\circ$,

 $h_p(a_\sigma)$ is a crossing point of $p$ and $q_1$ if
  $h_p(a_\sigma)\in \range(h_{q_1})$ \\
and correspondingly,

   $h_p(a_\sigma)$ is a crossing point of $p$ and $q_2$ if
$h_p(a_\sigma)\in \range(h_{q_2})$,

$h_p(b_\sigma)$ is a crossing point of $p$ and $q_1$ if
$h_p(b_\sigma)\in \range(h_{q_1})$ and

$h_p(b_\sigma)$ is a crossing point of $p$ and $q_2$ if
$h_p(b_\sigma)\in \range(h_{q_2})$.

\medskip
\noindent On the other hand, every crossing point of $p$ and $q_1$ or $q_2$ is equal to $h_p(a_\sigma)$ or $h_p(b_\sigma)$ for some $\sigma$ with $h_p(\sigma)\in\Delta^\circ$. Therefore, the number $N$ of
 crossings of $p$ with $q_1$ or $q_2$ is even. Since by (\ref{f30})
 $\{y,y'\}\cap \range(h_p)=\emptyset$,
$N={\rm CN}(p,q_1)+{\rm CN}(p,q_2)$ is an even number, hence
 $\pi(p,q_1)=\pi(p,q_2)$.

{\bf \boldmath Case (d) $y\neq y'$, $z_2\in  l(y,y')\setminus \overline {yy'}$ and $z_1\not\in l(y,y')$:} Let $q_3:=(r,y)(t_1,z_3)(r',y'))$ for
$z_3:=(z_1+z_2)/2$. By Case~(b), $\pi(p,q_2)=\pi(p,q_3)$ and by Case~(c), $\pi(p,q_1)=\pi(p,q_3)$.

{\bf \boldmath Case (e) $y\neq y'$,  $z_1,z_2\not\in l(y,y')$:}
Let $q_3:=(r,y)(t_1,z_3)(r',y'))$ for
$z_3:=(z_1+z_2)/2$. By Case~(c), $\pi(p,q_2)=\pi(p,q_3)$ and $\pi(p,q_1)=\pi(p,q_3)$.

{\bf \boldmath Case (f) $y\neq y'$ and $z_1,z_2\in l(y,y')\setminus\overline{yy'}$}: Proof via $q_3$ as in (d) and (e).
\qq\\

For tracks $p,q$ and $q'$, by tiny shifts of the vertices of $p$ we can obtain a track $\overline p$ such that
$\overline p \bowtie  q$ and $\overline p \bowtie  q'$.

\begin{lemma} \label{l4}$ $
\begin{enumerate}
\item  \label{l4a}
Let $p=((s_0,x_0),(s_1,x_1).\ldots,(s_k,x_k))$ and let
$q$ and $q'$ be tracks. Then for every $\delta>0$ there is some track $\overline p= ((s_0,y_0),(s_1,y_1),\ldots,(s_k,y_k)) $ such that
$\overline p \bowtie  q$, $\overline p \bowtie  q'$ and
$\|x_i-y_i\|<\delta$ for $0\leq i\leq k$.

\item  \label{l4b}
If $p=((s_0,x_0),(s_1,x_1).\ldots,(s_k,x_k))$  and
$\overline p= ((s_0,y_0),(s_1,y_1),\ldots,(s_k,y_k)) $ are tracks such that $\|x_i-y_i\|<\delta$ for $0\leq i\leq k$ then
$\|h_p(s)-h_{\overline p}(s)\| < \delta$ for all $s_0\leq s\leq s_k$.
\end{enumerate}
\end{lemma}

\proof \\
\ref{l4a}.\  We must find (prove the existence of) points $y_0,y_1,\ldots,y_k $ such that
$\|x_i-y_i\|<\delta$ and for $\overline p:=(s_0,y_0),(s_1,y_1),\ldots, (s_k,y_k))$,
\begin{eqnarray}
\mathcal V(\overline p)\cap (\mathcal L(q)\cup \mathcal L(q')) & = & \emptyset \an\\
\mathcal L(\overline p)\cap (\mathcal V(q)\cup \mathcal V(q')) & = & \emptyset
\end{eqnarray}
The following statements are true since $\mathcal V(q)\cup \mathcal V(q')$ is a finite set of points and $\mathcal L(q)\cup \mathcal L(q')$ is a finite set of straight lines.\\
-- There is some  $y_0\in B(x_0,\delta)$ such that $y_0\not\in \mathcal L(q)\cup \mathcal L(q')$.\\
-- Suppose $y_i$ has been determined for some $0\leq i<k$. There is some $y_{i+1}\in B(x_{i+1},\delta)$ such that
$$y_{i+1}\neq y_i,\ \ y_{i+1}\not\in (\mathcal L(q)\cup \mathcal L(q'))
\an  l(y_i,y_{i+1})\cap(\mathcal V(q)\cup\mathcal V(q'))=\emptyset ,\ $$

Then   $\overline p := ((s_0,y_0),(s_1,y_1).\ldots,(s_k,y_k))$ is a track such that $\overline p\bowtie q$, $\overline p\bowtie q'$ and
$\|x_i-y_i\|<\delta$ for $0\leq i\leq k$.\\

\ref{l4b}. \ For $s_i\leq s\leq s_{i+1}$ by (\ref{f5}),
\begin{eqnarray*}
\|h_p(s)-h_{\overline p}(s)\|
 & = & \|h_p(s_i)+\frac{s-s_i}{s_{i+1}-s_i}
 (h_p(s_{i+1})-h_p(s_i))\\
 &&-h_{\overline p}(s_i)-\frac{s-s_i}{s_{i+1}-s_i}
 (h_{\overline p}(s_{i+1})                     -h_{\overline p}(s_i))\|\\
&=& \|x_i+\frac{s-s_i}{s_{i+1}-s_i}
 (x_{i+1}-x_i)\\
 &&-y_i-\frac{s-s_i}{s_{i+1}-s_i}
 (y_{i+1}-y_i)\|\\
& = &  \|(1-\frac{s-s_i}{s_{i+1}-s_i})\cdot(x_i-y_i)\\
&& + \frac{s-s_i}{s_{i+1}-s_i})
\cdot(x_{i+1}-y_{i+1})\|\\
&< & \delta
\end{eqnarray*}
\qq \\

{\bf \boldmath In the following let $f,g:[-1;2]\to\IR^2$, let
$I=[a_I;b_I]$ and $J=[a_J;b_J]$ be intervals with rational endpoints such that $I,J\In [-1;2]$  and let $\alpha_{IJ}$ be the number from  (\ref{f8}).}

Since $f$ and $g$ are continuous, there is  a  modulus of uniform continuity ${\rm md}:\IN\to\IN$
for $f$ and $g$, that is, ${\rm md}$ is increasing and for all $t,t'\in[-1;2]$,
\begin{eqnarray}\label{f16}
\|f(t)-f(t')\|<2^{-n} \an \|g(t)-g(t')\|<2^{-n}&\mbox{if} &|t-t'|<2^{-{\rm md}(n)}.
\end{eqnarray}

We approximate $f|_I$ and $g|_J$ by tracks.

\begin{definition}\label{d4}
A track $p=((s_0,x_0),(s_1,x_1),\ldots, (s_k,x_k))$   is an
approximation with precision $2^{-n}$ (shortly an $n$-approximation) of $f|_I$  if
\begin{eqnarray}
\label{f7}&  s_0=a_I \an s_k=b_I\,,\\
\label{f3}&  s_{i+1}-s_i< 2^{-{\rm md}(n)}\an\\
\label{f4} & \|f(s_i)-x_i\|<2^{-n} \,.
\end{eqnarray}
\end{definition}
Remember that by Definition~\ref{d1}, $x_i=h_p(s_i)$. Approximations of $g|_J$ are defined accordingly.

\begin{lemma}\label{l1}
 In Definition~\ref{d4}, for all meaningful $i$ and $s$,
 \begin{eqnarray}
 \label{f9} \|x_i-x_{i+1}\| & < & 3\cdot 2^{-n}\,,\\
\label{f10} \|f(s)-h_p(s)\|&<&5\cdot 2^{-n}\,.
%\label{f11}d_I(f|_I,h_p)&<&5\cdot 2^{-n}\,.
  \end{eqnarray}
\end{lemma}
%where $d_I(f|_I,h_p):=\max_{ s\in I}\|f(s)-h_p(s)\|$.\\

\proof
 $ \|x_i-x_{i+1}\|\leq \|x_i-f(s_i)\| + \|f(s_i)-f(s_{i+1})\|
+\| f(s_{i+1})-x_{i+1})\|
< 3\cdot 2^{-n}$.

For $s_i\leq s\leq s_{i+1}$, \ \
$\|f(s)-h_p(s)\| \leq \|f(s)-f(s_i)\|+\|f(s_i)-h_p(s_i)\|+\|h_p(s_i)-h_p(s)\|
 < 2^{-n} + 2^{-n} +\|x_i-x_{i+1}\|<5\cdot 2^{-n} $
 by (\ref{f5}) and (\ref{f9}).
 \qq\\

After these technical preparations  prove the following central lemma.
\begin{lemma}\label{l6}
Let $2^{-n}<\alpha_{IJ}/16$.
 Let  $p$ and $p'$ be $n$-approximations of $f|_I$ and let $q$ and $q'$ be $n$-approximations of  $g|_J$ such that $p\bowtie q$ and $p'\bowtie q'$. Then $\pi(p,q)=\pi(p',q')$.
 \end{lemma}

\begin{figure}[htbp]
\setlength{\unitlength}{.24ex}
\begin{picture}(200,37)(-75,15)

\put(100,20){\line(2,1){40}}
\put(90,22){\line(0,1){16}}

\put(200,20){\line(-2,1){40}}
\put(210,22){\line(0,1){16}}

\put(90,15){\makebox(0,0)[cc]{$q$}}
\put(90,45){\makebox(0,0)[cc]{$p$}}
\put(150,45){\makebox(0,0)[cc]{$\overline p$}}
\put(210,45){\makebox(0,0)[cc]{$p'$}}
\put(210,15){\makebox(0,0)[cc]{$q'$}}

\end{picture}
\caption{The $\bowtie$-relation between the tracks $p,q,\overline p,p'$ and $q'$\,.} \label{fig6}
\end{figure}
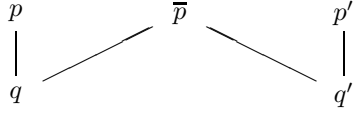

\proof
Let $$p=((s_0,x_0),(s_1,x_1).\ldots,(s_k,x_k))$$ where $s_0=a_I$ and $s_k=b_I$. Notice that possibly $p$ and $q'$ as well as $p'$ and $q$ have common straight line segments such that $\pi(p,q')$ and
$\pi(p',q)$  are not defined.
By Lemma~\ref{l4} there is some track
$$\overline p =((s_0,y_0),(s_1,y_1),\ldots,(s_k,y_k)) $$
such that
\begin{eqnarray}\label{f40}
\overline p \bowtie  q, & \overline p \bowtie  q' &\an
\|x_i-y_i\|< \ee \mbox{ for }  0\leq i\leq k\,.
\end{eqnarray}
Figure~\ref{fig6} shows the $\bowtie$-relation (indicated by lines) between the tracks $p,q,\overline p,p'$ and $q'$.
It suffices to prove
\begin{eqnarray}\label{f18}
\pi(p,q)=\pi(\overline p,q)=\pi(\overline p,q')
=\pi( p',q')\,.
\end{eqnarray}
We use the facts that by the definitions the points $f(t)$, $h_p(t)$, $h_{\overline p}(t)$ and $h_{p'}(t)$ are close together and that the points $g(t)$, $h_q(t)$ and $h_{q'}(t)$ are close together.

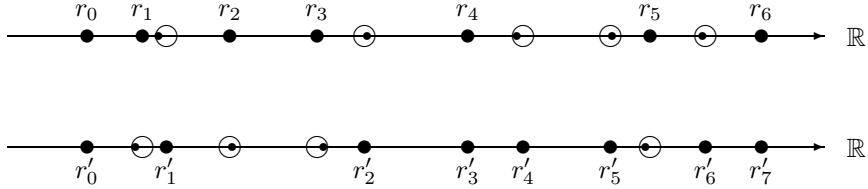
\begin{figure}[htbp]
\setlength{\unitlength}{3pt}
\begin{picture}(100,27)(-22,4)

\put(-3,10){\vector(1,0){103}}
\put(-3,24){\vector(1,0){103}}
\put(104,10){\makebox(0,0)[cc]{$\IR$}}
\put(104,24){\makebox(0,0)[cc]{$\IR$}}

\put(07,10){\circle*{1.7}}
\put(17,10){\circle*{1.7}}
\put(42,10){\circle*{1.7}}
\put(55,10){\circle*{1.7}}
\put(62,10){\circle*{1.7}}
\put(73,10){\circle*{1.7}}
\put(85,10){\circle*{1.7}}
\put(92,10){\circle*{1.7}}

\put(17,24){\circle{2.5}}  \put(16,24){\circle*{1}}
\put(42,24){\circle{2.5}}  \put(42.3,24){\circle*{1}}
\put(62,24){\circle{2.5}}  \put(61.2,24){\circle*{1}}
\put(73,24){\circle{2.5}}  \put(73.2,24){\circle*{1}}
\put(85,24){\circle{2.5}}  \put(84.6,24){\circle*{1}}

\put(07,7){\makebox(0,0)[cc]{$r'_0$}}
\put(17,7){\makebox(0,0)[cc]{$r'_1$}}
\put(42,7){\makebox(0,0)[cc]{$r'_2$}}
\put(55,7){\makebox(0,0)[cc]{$r'_3$}}
\put(62,7){\makebox(0,0)[cc]{$r'_4$}}
\put(73,7){\makebox(0,0)[cc]{$r'_5$}}
\put(85,7){\makebox(0,0)[cc]{$r'_6$}}
\put(92,7){\makebox(0,0)[cc]{$r'_7$}}

\put(07,24){\circle*{1.7}}
\put(14,24){\circle*{1.7}}
\put(25,24){\circle*{1.7}}
\put(36,24){\circle*{1.7}}
\put(55,24){\circle*{1.7}}
\put(78,24){\circle*{1.7}}
\put(92,24){\circle*{1.7}}

\put(14,10){\circle{2.5}}
\put(25,10){\circle{2.5}}
\put(36,10){\circle{2.5}}
\put(78,10){\circle{2.5}}

\put(13.1,10){\circle*{1}}
\put(25.3,10){\circle*{1}}
\put(36.8,10){\circle*{1}}
\put(77.4,10){\circle*{1}}

\put(07,27){\makebox(0,0)[cc]{$r_0$}}
\put(14,27){\makebox(0,0)[cc]{$r_1$}}
\put(25,27){\makebox(0,0)[cc]{$r_2$}}
\put(36,27){\makebox(0,0)[cc]{$r_3$}}
\put(55,27){\makebox(0,0)[cc]{$r_4$}}
\put(78,27){\makebox(0,0)[cc]{$r_5$}}
\put(92,27){\makebox(0,0)[cc]{$r_6$}}

\end{picture}
\caption{Insertion of redundant vertices into $q$ and $q'$.} \label{fig3}
\end{figure}

In the following we prove the second equation $\pi(\overline p,q)=\pi(\overline p,q')$ of (\ref{f18}) in detail.
Let
\begin{eqnarray}
\label{f19}  q & = &((r_0,u_0),(r_1,u_1),\ldots,(r_{l-1},u_{l-1}), (r_l,u_l))\,,\\
\label{f20}  q' & = &((r'_0,u'_0),(r'_1,u'_1),\ldots,(r'_{m-1},u'_{m-1}), (r'_m,u'_m))\ \\
\label{f21} &\ \mbox{where} & r_0=r'_0=a_J \an  r_l=r'_m=b_J\,.
\end{eqnarray}
be $n$-approximations of $g|_J$.  As an example, Figure~\ref{fig3} shows the arguments $r_0\ddd r_6$  of a track $q$ on the real line and and the arguments $r'_0\ddd r'_7$ of a track $q'$ on the real line (thick black dots).

From (\ref{f40}) we know that $\mathcal V(q)\cap \lp=\emptyset$ and $\mathcal V(q')\cap \lp=\emptyset$.
In a first step we add ``redundant" vertices to $q$ and $q'$ such that the resulting tracks $q_1$ and $q'_1$ have the same number of vertices and such that $h_q=h_{q_1}$, $h_{q'} = h_{q'_1}$,
$\mathcal V(q_1)\cap \lp=\emptyset$ and $\mathcal V(q'_1)\cap \lp=\emptyset$.

Let $i$ be a number such that $r_i\not\in \{r'_0\ddd r'_m\}$. There is a unique $j$ such that $r'_j<r_i<r'_{j+1}$ (e.g. $i=3$ and $j=1$ in Figure~\ref{fig3}).\\

The point $h_{q'}(r_i)$ is an element of the straight line segment
$\overline{u'_ju'_{j+1}}$ and different from $u'_j$ and $u'_{j+1}$.
We would like to add the ``redundant" pair $(r_i,h_{q'}(r_i))$ to $q'$ with result
$$((r'_0,u'_0)\ddd (r'_j,u'_j),(r_i,h_{q'}(r_i)(r'_{j+1},u'_{j+1})\ddd (r'_m,u'_m))$$
But possibly $h_{q'}(r_i)\in \lp$. Instead, we add a pair $(t,v)$ with $v=h_{q'}(t)\not\in  \lp$ to $q'$ which is very close to $(r_i,h_{q'}(r_i))$.
Let
\begin{eqnarray}\label{f43}
\gamma:= \min\{|c-d|\mid  c,d\in\{
r_0\ddd r_l,r'_0\ddd r'_m\} \an c\neq d \}\,.
\end{eqnarray}
Since $u'_j\not\in \lp$, by Lemma~\ref{l7}
 $\: l(u'_ju'_{j+1})$  intersects every straight line from
 $\lp$ at most once. Therefore, the straight line segment $\overline {u'_j u'_{j+1}}$ contains only finitely many points of
$\lp$.
Thus there is some $t$ such that $|t-r_i|<\gamma/8$ and $v:=h_{q'}(t)\not\in\lp$.
Then $v\in \overline{u'_ju'_{j+1}}$ (see (\ref{f5})).
We add the pair $(t,v)$ to $q'$ with result
$$\overline q:=((r'_0,u'_0)\ddd (r'_j,u'_j),(t,h_{q'}(t)(r'_{j+1},u'_{j+1})\ddd (r'_m,u'_m))\,.$$
Then
$h_{q'}=h_{\overline q}$ .
Let $q'_1$ be the track obtained from $q'$ by adding a pair
$(t,v) $  in this way for every $i$ such that  $r_i\not\in \{r'_0\ddd r'_m\}$ in turn. Then $h_{q'}=h_{q'_1}$ and $\mathcal V(q'_1)\cap \lp=\emptyset$.

Correspondingly let  $q_1$ be the track obtained from $q$
in the same way. Then $h_{q}=h_{q_1}$ and
 $\mathcal V(q_1)\cap \lp=\emptyset$. In summary,
\begin{eqnarray}\label{f34}
h_{q}=h_{q_1}, \  h_{q'}=h_{q'_1}, \ \ \mathcal V(q_1)\cap \lp=\emptyset \an \mathcal V(q'_1)\cap \lp=\emptyset\,.
\end{eqnarray}

In Figure~\ref{fig3} the inserted arguments $t$  are within the  the circles.
The new tracks can be written as
\begin{eqnarray}
\label{f41} q_1 & = & ((t_0, v_0),(t_1,v_1)\ddd (t_\mu, v_\mu))\,,\\
\label{f42} q'_1 & = & ((t'_0, v'_0),(t'_1,v'_1)\ddd (t'_\mu, v'_\mu))
\end{eqnarray}
where $a_J=r_0=t_0=t'_0$ and $b_J=r_l=t_\mu=t'_\mu$.
By the condition $|t-r_i|<\gamma/8$ in the definition of $(t,v)$ above,
%
%\begin{eqnarray}\label{f47}
%&|t_\nu-t'_\nu|<\gamma/8,\ \ t_\nu<t_{\nu+1},\ \  \ t'_\nu<t_{\nu+1},\ \
% \ t_\nu<t'_{\nu+1},\ \  \ t'_\nu<t'_{\nu+1}\,.
%\end{eqnarray}
%
\begin{eqnarray}\label{f47}
&|t_\nu-t'_\nu|<\gamma/8 \an \{t_\nu,t_\nu'\}< \{t_{\nu+1},t'_{\nu+1}\}\,. \end{eqnarray}
The next Proposition prepares the proof or Proposition~\ref{prop3}.

\begin{proposition}\label{prop2}
For the tracks $q_1$ and $q'_1$, for all $ 0<\nu<\mu$,
\begin{eqnarray}
\label{f35} &\{v_{\nu-1},v_\nu,v_{\nu+1}, v'_{\nu-1},v'_\nu,v'_{\nu+1}\}\In B(g(t_\nu),5\eee) \an \\
\label{f54} & \{y_0,y_k\}\cap B(g(t_\nu),14\eee)=\emptyset\,,\\
\label{f84} & B(g(t_1),9\eee)\cap \range(h_{\overline p})=\emptyset\,,\\
\label{f85} & B(g(t_{\mu-1}),9\eee)\cap \range(h_{\overline p})=\emptyset\,.
\end{eqnarray}
\end{proposition}

\proof (Proposition~\ref{prop2})

 {\bf Proof of (\ref{f35}):}
 Consider $((t_{\nu-1},v_{\nu-1}), (t_{\nu},v_{\nu}),(t_{\nu+1},v_{\nu+1}))$ as a part of $q_1$. Then $(t_{\nu},v_{\nu})$ has already been in $q$,
 that is,  $(t_\nu,v_\nu)=(r_i,u_i)$ for some $i$ or it has been inserted via  some~$j$ such that $r'_j\not\in \{r_0\ddd r_l\}$. Therefore, there is some $i$ such that (see Figure~\ref{fig3})
\begin{eqnarray}
\label{f49}&r_{i-1}\leq t_{i-1} <r_i< t_{\nu+1}\leq r_{i+1}\ \ \mbox{or}\\
\label{f50}&r_i\leq t_{\nu-1}<t_\nu<t_{\nu+1}\leq r_{i+1}\,.
\end{eqnarray}

Consider (\ref{f49}). Then
$\|g(t_\nu) -v_\nu\| =\|g(r_i)-u_i\| < \ee$ by (\ref{f4}). Furthermore,
$\|g(t_\nu) -v_{\nu-1}\|\leq \|g(t_\nu) -v_\nu\| + \|v_\nu -v_{\nu-1}\|
<\ee + \|u_i- v_{\nu-1}\|$. Since $h_q=h_{q_1}$,
$v_{\nu-1}\in \overline{u_{i-1}u_i}$, hence
$\|u_i- v_{\nu-1}\|\leq \|u_{i-1}-u_i\| <3\eee $ by (\ref{f9}).
Therefore, $\|g(t_\nu) -v_{\nu-1}\|<4\eee$. By symmetry
$\|g(t_\nu) -v_{\nu+1}\|<4\eee$.
Therefore,
$\{v_{\nu-1},v_\nu,v_{\nu+1}\}\In B(g(t_\nu),4\eee)$.

Consider (\ref{f50}).Then

$\|g(t_\nu)-v_\nu\| = \|g(t_\nu)- h_q(t_\nu)\| <5\eee $ by(\ref{f10}).
Furthermore,
$\|g(t_\nu)-v_{\nu-1}\| \leq  \|g(t_\nu)-v_\nu\| + \|v_\nu- v_{\nu-1}\|
< \ee + \| u_{i+1}-u_i\| < \ee + 3 \eee =4\eee $
(since $v_\nu, v_{\nu-1}\in \overline{u_iu_{i+1}}$ and by (\ref{f9})).
By symmetry
$\|g(t_\nu) -v_{\nu+1}\|<4\eee$.
Therefore,
$\{v_{\nu-1},v_\nu,v_{\nu+1}\}\In B(g(t_\nu),4\eee)$.

As a summary, in both cases (\ref{f49}) and (\ref{f50}),
$\{v_{\nu-1},v_\nu,v_{\nu+1}\}\In B(g(t_\nu),4\eee)$.
By symmetry, $\{ v'_{\nu-1},v'_\nu,v'_{\nu+1}\}\In B(g(t'_\nu),4\eee)$.
Since $|t_\nu-t'_\nu|<\gamma/4 <2^{-{\rm md}(n)}$ , $\|g(t_\nu)-g(t'_\nu\| < \ee$ by (\ref{f47}), (\ref{f43}), (\ref{f3}) and (\ref{f16}).
Therefore, $\{v_{\nu-1},v_\nu,v_{\nu+1}, v'_{\nu-1},v'_\nu,v'_{\nu+1}\}\In B(g(t_\nu),5\eee)$.

{\bf Proof of (\ref{f54}):} By (\ref{f8}), (\ref{f4}) and
Lemma~\ref{l4}  for all $\nu$

$16\eee<\alpha_{IJ} \leq \|f(s_0)-g(t_\nu)\|
\leq \|f(s_0)-x_0\|+\|x_0-y_0\| +\|y_0-g(t_\nu)\|
<\ee +\ee +\|y_0-g(t_\nu)\|$, hence  $\|y_0-g(t_\nu)\| > 14 \eee$. By symmetry,
 $\|y_k-g(t_\nu)\| > 14 \eee$.

 {\bf Proof of (\ref{f84}):}
 For every $s_0\leq s\leq s_k$, $\|f(s)- g(r_0)\|\geq\alpha_{IJ}>16\eee$. By (\ref{f10}) and Lemma~\ref{l4},
  $\|h_{\overline p}(s)- g(r_0)\|>10\eee$. Since $|r_0-t_1|\leq |s_0-s_1|<2^{\rm md}(n)$, $\|h_{\overline p}(s)- g(t_1)\|>9\eee$,
hence $B(g(t_1),9\eee)\cap \range(h_{\overline p})=\emptyset$.

 {\bf Proof of (\ref{f85}):} From (\ref{f84}) by symmetry.
\phantom{$ $} \hfill $\Box$~(Proposition~\ref{prop2})
\medskip

\begin{proposition}\label{prop3} \ \
$\pi(\overline p,q_1)=\pi(\overline p,q'_1)\,.$
\end{proposition}

\proof (Proposition~\ref{prop3}):  We transform the track $q'_1$  in two phases to the track $q_1$ preserving the crossing parity in each step.
Remember (\ref{f34}) and Proposition~\ref{prop2} for
\begin{eqnarray}
\nonumber q_1 & = & ((t_0, v_0),(t_1,v_1),(t_2,v_2)\ddd (t_\mu, v_\mu)) \an\\
\nonumber q'_1 & = & ((t'_0, v'_0),(t'_1,v'_1),(t'_2,v'_2)\ddd (t'_\mu, v'_\mu))\,.
\end{eqnarray}

In the first phase for every $0\leq i<\mu$ we replace every $v'_i$ by some $w_i$ without changing $q'_i$.
 For $0\leq i < \mu$ let $Q(i)$ be the following property:\\
 There are points $w_0\ddd w_i$ such that for all $ 1\leq j\leq i$
\begin{eqnarray}
%\label{f80}&& w_i\neq v'_{i+1}, \\
\label{f81}& %\forall( 1\leq j\leq i)(
\|w_j-v'_j\|< \ee/16,\ \ w_j \not\in\lp, \ w_j\not\in\{v_{j-1}, w_{j-1}, v'_{j+1}\}) \an \\
\label{f82}&\ \  \pi(\overline p,q'_1)= \pi(\overline p,\overline q_i)
\end{eqnarray}
where $\overline q_i$ is the track
$$\overline q_i:=((t'_0, w_0)\ddd (t'_i,w_i),(t'_{i+1},v'_{i+1})\ddd (t'_\mu, v'_\mu))\,.$$

Let $w_0:= v'_0$. Then $q'_1=\overline q_0$, hence
$Q(0)$ is true.
Suppose for some $0\leq i \leq \mu-2$ we have proved $Q(i)$.
There is some $w_{i+1}$ such that
\begin{eqnarray}\label{f79}
\|w_{i+1}-v'_{i+1}\|< \ee/16,\ \ w_{i+1}\not\in\lp,\ \  w_{i+1}\not\in\{v_i, w_i,v'_{i+2}\}\,.
\end{eqnarray}
By definition,
\begin{eqnarray*}
\overline q_i & = & ((t'_0, w_0)\ddd \mbox{\boldmath$(t'_i,w_i),(t'_{i+1},v'_{i+1})
,(t'_{i+2},v'_{i+2})$}\ddd (t'_\mu, v'_\mu))\an\\
\overline q_{i+1} &=& ((t'_0, w_0)\ddd \mbox{\boldmath $(t'_i,w_i),(t'_{i+1},w_{i+1}),
(t'_{i+2},v'_{i+2})$}\ddd (t'_\mu, v'_\mu))\,.
\end{eqnarray*}
Since $\overline q_i$ is a track, by (\ref{f79}), $\overline q_{i+1}$
is a track (neighboring vertices must be different).
Also  by (\ref{f79}), (\ref{f81}) is true also for $i+1$.
By Proposition~\ref{prop2} Lemma~\ref{l8} can be applied to the boldface sub-tracks $\tau_1:=((t'_i,w_i),(t'_{i+1},v'_{i+1})
,(t'_{i+2},v'_{i+2}))$ and $\tau_2:=((t'_i,w_i),(t'_{i+1},w_{i+1}),
(t'_{i+2},v'_{i+2}))$ such that $\pi(\overline p,\tau_1)=\pi(\overline p,\tau_2)$. Since the vertices of  $\overline q_i$ and   $\overline q_{i+1}$ are not in $\mathcal L(\overline p)$, hence not in $\range(h_{\overline p})$,
$\pi(\overline p, q'_1)=\pi(\overline p,\overline q_i)
=\pi(\overline p,\overline q_{i+1})$. Therefore, we have proved $Q(i+1)$.

By induction $Q(\mu-1 )$ is true, hence
\begin{eqnarray}\label{f83}
\pi(\overline p, q'_1)
=\pi(\overline p,\overline q_{\mu-1}).
\end{eqnarray}

\begin{figure}[htbp]
\setlength{\unitlength}{2.45pt}
\begin{picture}(140,25)(-26,-2)

\multiput(0,0)(20,0){7}{\line(1,0){10}}
\multiput(0,20)(20,0){7}{\line(1,0){10}}
\multiput(0,15)(20,0){7}{\line(1,-1){10}}
\multiput(70,10)(2.5,0){5}{\circle*{1}}
\put(0,0){\line (1,0){10}}
\put(20,0){\line (1,0){10}}
\put(120,0){\line (1,0){10}}

\put(-3,0){\makebox(0,0)[cc]{$w_0$}}
\put(15,0){\makebox(0,0)[cc]{$w_1$}}
\put(35,0){\makebox(0,0)[cc]{$w_2$}}
\put(55,0){\makebox(0,0)[cc]{$w_3$}}
\put(95,0){\makebox(0,0)[cc]{$w_{\mu-2}$}}
\put(115,0){\makebox(0,0)[cc]{$w_{\mu-1}$}}
\put(135,0){\makebox(0,0)[cc]{$v'_\mu$}}

\put(-3,20){\makebox(0,0)[cc]{$v_0$}}
\put(15,20){\makebox(0,0)[cc]{$v_1$}}
\put(35,20){\makebox(0,0)[cc]{$v_2$}}
\put(55,20){\makebox(0,0)[cc]{$v_3$}}
\put(95,20){\makebox(0,0)[cc]{$v_{\mu-2}$}}
\put(115,20){\makebox(0,0)[cc]{$v_{\mu-1}$}}
\put(135,20){\makebox(0,0)[cc]{$v_\mu$}}

\end{picture}
\caption{ The $\neq$-relation for the vertices of
$q_1$ (top) and
$\overline q_{\mu-1}$ (bottom).} \label{fig9}
\end{figure}
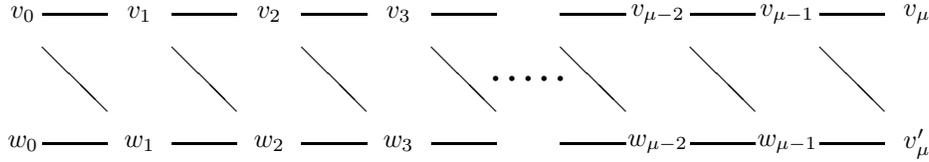

Figure~\ref{fig9} shows the $\neq$-relation for the vertices of\\
$q_1=((t_0, v_0),(t_1,v_1),(t_2,v_2)\ddd (t_\mu, v_\mu))$ (top) and\\
$\overline q_{\mu-1}= ((t'_0, w_0),(t'_1,w_1)\ddd (t'_{\mu-1},w_{\mu-1}),(t'_\mu, v'_\mu))$
(bottom). A line between $x$ and $y$ means $x\neq y$.

In the second phase two we transform $\overline q_{\mu-1}$ step by step to $q_1$ without changing the crossing parity with $\overline p$.
For $0\leq i\leq \mu-1$ let
$$q''_i:=((t_0, v_0)\ddd (t_i,v_i)(t'_{i+1},w_{i+1})\ddd
(t'_{\mu-1},w_{\mu-1}),(t'_\mu, v'_\mu))\,.$$
Since $q_1$ and $\overline q_{\mu-1}$ are tracks and $v_i\neq w_{i+1}$,
$q''_i$ is a track.

By (\ref{f35})  and by $\|w_1-v'_1\|<\ee/16$ (\ref{f81}),
$\overline{ v_0w_1}\cup \overline {w_0w_1}\In B(g(t_1),6\eee)$.
By (\ref{f84}),
$\overline{ v_0w_1}$ and $\overline {w_0w_1}$ do not intersect
$h_{\overline p}\,$. Therefore,
\begin{eqnarray}\label{f86}
\pi(\overline p, \overline q_{\mu-1})&=& \pi(\overline p,q''_0)
\end{eqnarray}

Suppose for $0\leq i\leq \mu-3$
$$\pi(\overline p, \overline q_{\mu-1}) = \pi(\overline p,q''_i)\,.$$
By definition,
\begin{eqnarray*}
q''_i &=& ((t_0, v_0)\ddd \mbox{\boldmath $(t_i,v_i),(t'_{i+1},w_{i+1}),(t'_{i+1},w_{i+2})$}\ddd
(t'_\mu, v'_\mu)) \an \\
q''_{i+1}&=& ((t_0, v_0)\ddd \mbox{\boldmath$(t_i,v_i),(t_{i+1},v_{i+1}),(t'_{i+1},w_{i+2})$}\ddd
(t'_\mu, v'_\mu)) \,,
\end{eqnarray*}

By Proposition~\ref{prop2} Lemma~\ref{l8} can be applied to the boldface sub-tracks

$\tau_1:=((t_i,v_i),(t'_{i+1},w_{i+1}),(t'_{i+1},w_{i+2})$ and

$\tau_2:=((t_i,v_i)(t_{i+1},v_{i+1}),(t'_{i+1},w_{i+2}))$

such that $\pi(\overline p,\tau_1)=\pi(\overline p,\tau_2)$.

 Since the vertices of  $\overline q''_i$ and   $\overline q''_{i+1}$ are not in $\mathcal L(\overline p)$, hence not in $\range(h_{\overline p})$,
$\pi(\overline p, \overline q_{\mu-1})=\pi(\overline p, q''_i)
=\pi(\overline p,q''_{i+1})$. By induction,

$$\pi(\overline p, \overline q_{\mu-1})=\pi(\overline p, q''_{\mu-2})$$
Finally we compare
\begin{eqnarray*}
q_1 & = & ((t_0,v_0),\ddd (t_{\mu-2},v_{\mu-2}),(t_{\mu-1},v_{\mu-1}),(t_\mu,v_\mu))\an\\
q''_{\mu-2} & = & ((t_0,v_0),\ddd (t_{\mu-2},v_{\mu-2}),(t'_{\mu-1},w_{\mu-1}),(t'_\mu,v'_\mu))\,.
\end{eqnarray*}

Since $\|w_{\mu-1}-v'_{\mu-1}\|<\ee$ by Proposition~\ref{prop2} and
$$(\overline{v_{\mu-2}v_{\mu-1}}\cup \overline{v_{\mu-1}v_\mu}
\cup \overline{v_{\mu-2}w_{\mu-1}}\cup \overline{w_{\mu-1}v'_\mu})
\cap \range(h_{\overline p})=\emptyset\,,$$

$\pi(\overline p, q''_{\mu-2})=\pi(\overline p,q_1)$, hence
$\pi(\overline p, q'_1)=\pi(\overline p,q_1)$.
\phantom{$ $} \hfill $\Box$~(Proposition~\ref{prop3})
\bigskip

\noindent We continue the proof of Lemma~\ref{l6}. Since $h_q=h_{q_1}$ and $h_{q'}=h_{q'_1}$ by (\ref{f34}),
$\pi(\overline p, q')=\pi(\overline p,q)$, which is the second equation from (\ref{f18}). The first and the third equation can be poved accordingly. We omit the details.
\phantom{$ $} \hfill $\Box$~(Lemma~\ref{l6})
\\

Lemma~\ref{l6} allows to define the parity of the pair $(f|_I,g|_J)$.

\begin{definition}\label{d2}
$\pi(f|_I,g|_J):=\pi(p,q)$ for an arbitrary $n$-approximation $p$ of $f|_I$ and an arbitrary $n$-approximation $q$ of $g|_J$ such that $p \bowtie q$ and $\ee<\alpha_{IJ}/16$.
\end{definition}

\begin{lemma}\label{l5}
Let $f,g,J$ and $I=[a_I;b_I]$ be as before such that (\ref{f8}). Let $I_1:=[a_I;c]$ and $I_2:=[c;b_I]$ where $a_I<c<b_I$ and $f(c)\not\in g(J)$. Then
\begin{eqnarray}\label{f11}
\pi(f|_I,g|_J) =( \pi(f|_{I_1},g|_J) + \pi(f|_{I_2},g|_J))\;\mod\;2 \,.
\end{eqnarray}
\end{lemma}

\proof There are numbers $\alpha>0$ and $n\in \IN$ such that\\
$\alpha <d_s(\{f(a_I),f(c),f(b_I)\},g(J))$,
$\alpha <d_s(\{g(a_J),g(b_J)\},f(I))$ and
$16\eee <\alpha$.
There are $n$-approximations $p=((s_0,x_0)\ddd(s_k,x_x))$
and $q$ of $f|_I$ and $g|_q$, respectively, such that
 $c=s_i$ for some $i$ and $p\bowtie q$.

Let $p_1:=((s_0,x_0)\ddd(s_i,x_i))$ and
 $p_2:=((s_i,x_i)\ddd(s_k,x_k))$.

Then $p_1\bowtie q$ and $p_2\bowtie q$. Therefore
$ {\rm CN}(p_1,q)$ and ${\rm CN}(p_2,q)$ are well-defined such that
 (see Definition~\ref{d3})
 $${\rm CN}(p,q)= {\rm CN}(p_1,q) + {\rm CN}(p_2,q)\,,$$
hence
$\pi(p,q)={\rm CN}(p,q)\;\mod\; 2 =
( {\rm CN}(p_1,q)\;\mod\; 2 + {\rm CN}(p_2,q)\;\mod\; 2  )\;\mod\;2
=(\pi(p_1,q) + \pi(p_2,q) )\;\mod\;2$.
(\ref{f11}) follows by Definition~\ref{d2}.
\qq

\section {Crossing parities can be computed.}\label{secc}

In TTE
%for a sufficiently large finite set
%$\Sigma$ (the alphabet)
%computability on the set $\s$ of finite strings and on the set $\om$ of infinite sequences is defined explicitly by Turing %machines with finite or infinite one-way input and output tapes.
computability is defined on represented sets $(X,\delta)$.
A representation is a function $\delta:\In \s \to X$ (if $X$ is finite) or $\delta:\In \om \to X$ (if $X$ has at most continuum cardinality). If $\delta(w)=x$ then $w$ is considered as a ``name" (or a ``$\delta$-name") of $x$. Every $x\in X$ must have a name (and may have many names) but not every $w\in\s$ must be a name of some $x\in X$, hence $\delta$  is a partial surjective function.

For represented sets $(X_i,\delta_i)$ where $\delta_i: \In A_i\to X_i$, $A_i\in \{\s,\om\}$, a partial function $f:\In X_1\to X_2$ is computable, if there is a computable function $h:\In A_1\to A_2$
which realizes $f$, that is, if $w$ is a $\delta_1$-name of $x\in\dom(f)$ then $h(w)$ is a $\delta_2$-name of $f(x)$ (accordingly functions on Cartesian products).
For most sets in Analysis there are canonical (or standard or effective or obvious) representations which we use here.
We will say ``computable" without mentioning the (standard) representations.

We need the concept of computable multi-functions on represented sets. As an example, for the standard represented sets $(\IR,\rho)$, $(\IQ,\nu_\IQ)$ and $(\IN,\nu_\IN)$ there is no computable function $f:\IR\times \IN\to \IQ$ such that $|x-f(x,n)|< 2^{-n}$. But there is a computable function $h:\om\times \s\to\s$ which from every name of $x$ and every name of $n$ computes a name of some $a\in \IQ$ such that $|x-a|<2^{-n}$. The computed number $a$ depends on the names of $x$ and $n$ and not only on $(x,n)$.
The function $h$ realizes the
{\em multi-function} $f:\IR\times \IN\mto\IQ$ such that $f(x,n)=\{ a\in \IQ\mid |x-q|<2^{-n}\}$. Informally we say: $f$ maps every $(x,n)$ to some $a\in\IQ$ such that $|x-a|<2^{-n}$, more formally,  $f:(x,n)\mmto a$ such that $|x-a|<2^{-n}$.
In the following we apply results about computability on sets with standard representations from \cite{Wei00,BHW08,WG09,TW11}.

In Lemma~\ref{l3} we consider $n$-approximations $p$ and $q$ of $f|_I$ and $g|_J$, respectively such that $p\bowtie q$. We prove that $p$ and $q$ can be computed from $(f,g,I,J,n)$. We use from TTE, that $(f,s,n)\mmto y$ ($s\in\IQ, y\in\IQ^2$) such that $\|f(s)-y\| <\ee$) is computable.

\begin{lemma}\label{l3}
Let $\mathcal R:(f,g,I,J,n)\mmto (p,q)$ be the multi-function
mapping continuous functions $f,g:[-1;2]\to \IR^2$, closed rational intervals $I,J \In [-1;2]$ and a number $n\in\IN$ to a pair
$(p,q)$,  $ p = ((s_0,x_0)\ddd(s_k,x_k))$ and
$ q  = ((t_0,y_0)\ddd(t_l,y_l))$,
of rational tracks such that
\begin{eqnarray}
\label{f22}& \mbox{$p$ is an $n$-approximation of $f|_I$}\,,\\
\label{f23}& \mbox{$q$ is an $n$-approximation of $g|_J$} \an \\
%$\label{f24}& {\rm card}(\mathcal V(p)\cup \mathcal V(q))=k+l+2 \an\\
%\label{f24}& {\rm card}(\mathcal V(p))=k+1 \an  {\rm card}( \mathcal V(q))=l+1 \an\\
\label{f25}&  p\bowtie q\,.
\end{eqnarray}
 Then $\mathcal R$ is computable.
\end{lemma}
 Since rational tracks are ``discrete" objects there is no compuble funtion but only a computable multi-function (known from TTE).

\medskip
\proof
It is known \cite{Wei00,BHW08,WG09} that from continuous functions $f,g:[-1;2]\to\IR^2$ we can compute
a modulus of continuity ${\rm md}_f$ of $f$ and
a modulus of continuity ${\rm md}_g$ of $g$.
We can compute ${\rm md} :=\max ({\rm md}_f, {\rm md}_g)$, which is a modulus of continuity of $f$ and of $g$. Let  $I=[a_I;b_I]$ and $J=[a_J;b_J]$. Consider Definition~\ref{d4}.

First we compute~$p$.\\
Choose some $k$ and rational numbers $s_0:=a_I,\ s_1\ddd s_{k-1},s_k:=b_I$ such that $0<s_{i+1}-s_i<2^{-{\rm md}(n)}$ for
$0\leq i<k$.
Choose some $x_0\in\IQ^2$ such that $\|f(s_0)-x_0\|<\ee$.
Suppose, $x_i$ has been chosen for $i<k$. Then choose
some $x_{i+1}\in\IQ^2$ such that $\|f(s_{i+1})-x_{i+1}\|<\ee$ and $x_{i+1}\neq  x_i$.
Since $x_i\neq x_{i+1}$ for $i<k$,
$p:=((s_0,x_0)\ddd (s_k,x_k))$ is a rational track.
 It is an $n$-approximation of $f|_I$.

Next we compute $q$.\\
Choose some $l$ and rational numbers $t_0:=a_J,\ t_1\ddd t_{l-1},t_l:=b_J$ such that $0<t_{i+1}-t_i<2^{-{\rm md}(n)}$ for
$0\leq i<l$.
Choose some $y_0\in\IQ^2$ such that $\|g(t_0)-y_0\|<\ee$ and
$y_0\not\in  \mathcal L(p)$.
Suppose, $y_0\ddd y_i$ have been chosen for some $i<l$.
 Then choose
some $y_{i+1}\in\IQ^2$ such that
 $y_{i+1}\neq y_i$,
$\|g(t_{i+1})-y_{i+1}\|<\ee$,
$y_{i+1}\not\in  \mathcal L(p)\}$  and
$l(y_i,y_{i+1}))\cap \mathcal V (p)=\emptyset$.

Since $y_i\neq y_{i+1}$ for $0i < l$, $ q  = ((t_0,y_0)\ddd(t_l,y_l))$ is a rational track. It is an $n$-approximation of $g|_J$.
We have $\mathcal V(q)\cap \mathcal L(p)=\emptyset$ since $y_i\not\in \mathcal L(p)$ for all $i$. We have
$\mathcal L(q)\cap \mathcal V(p)=\emptyset$ since
 $l(y_i,y_{i+1}))\cap \mathcal V (p)=\emptyset$ for all $i<l$. Therefore, $p\bowtie q$.

All of this can be computed from
continuous functions $f,g:[-1;2]\to \IR^2$, rational intervals $I,J \In [-1;2]$ and $n$.
\qq

\begin{corollary}\label{cor1}
Let $f,g:[-1;2]\to \IR^2$ be continuous functions and let $I,J\In [-1;2]$ be closed rational intervals such that $f(I)\cap g(J)=\emptyset$.
Then $\pi(f|_I,g|_J)=0$.
\end{corollary}

\proof Since $f(I)$ and $ g(J)$ are compact, $d_s(f(I),g(J))>0$ and hence $\alpha_{IJ}>0$. There is some number $n$ such that $11\eee< d_s(f(I),g(J))$ and
 $\ee< \alpha_{IJ}/16$.

 By Lemma~\ref{l3} there are $n$-approximations $p$ and $q$ of $(f|_I$ and $g|_J$, respectively, such that $p \bowtie p$. By definition~\ref{d2}, $\pi(f|_I,g|_J)=\pi(p,q)$.
By (\ref{f10}), $\range(h_p)\In U(\range(f_I),5\eee)$ and
 $\range(h_q)\In U(\range(g_J),5\eee)$
\footnote{ where
 $U(A,\delta):=\{ x\in\IR^2\mid (\exists \,y\in A)\,\|x-y\|<\delta\}$}.
Since $11\eee< d_s(f(I),g(J))$,
$ U(\range(f_I),5\eee)\cap U(\range(g_J),5\eee)=\emptyset$, hence
$\range(h_p)\cap\range(h_q)=\emptyset$. We obtain ${\rm CN}(h_p,h_q)=0$
and $\pi(p,q)=0$. By Definition~\ref{d2}, $\pi(f|_I,g|_J)=0$.
\qq

\begin{lemma}\label{l13} $ $
\begin{enumerate}
\item \label{l13a}From continuous functions $f,g:[-1;2]\to \IR^2$  and  closed rational intervals $I,J\In [-1;2]$ we can compute $\alpha_{IJ}$.

\item \label{l13b} From continuous functions $f,g:[-1;2]\to \IR^2$  and  closed rational intervals $I,J\In [-1;2]$ such that $\alpha_{IJ}>0$ we
     can compute $\pi(f|_I,g|_J)$.
\end{enumerate}
\end{lemma}
 Remember that   $f(a_I),f(b_I)\not\in g(J) \an  g(a_J),g(b_J)\not \in f(I)$ is equivalent to $\alpha_{IJ}>0$.\\

\proof
(\ref{l13a}) From $f$, $a_I$ and $b_I$ we can compute the compact set $\{a_I,b_I\}$. from $g$ and $J$ we can compute the compact set $g(J)$. The function
$(A,B)\mapsto d_s(A,B)$ on (non-empty) compact set is computable.
Therefore, from $f,g,I,J$ we can compute $d_s(\{f(a_I),f(b_I)\}, g(J))$.
Correspondingly we can compute $d_s(\{g(a_J),g(b_J)\}, f(I))$ and hence their minimum $\alpha_{IJ}$.

(\ref{l13b}) Compute $\alpha_{IJ}$.
Find some $n$ such that $\ee < \alpha/16$. By Lemma~\ref{l3} we can compute $n$-approximations $p$ and $q$ of $f_I$ and $g_J$, respectively such that $p\bowtie q$. Compute the number of crossings ${\rm CN}(p,q)$ and its parity $\pi(p,q)$. By Definition~\ref{d2},
$\pi(f|_I,g|_J)=\pi(p,q)$
\qq

\section {The crossing parity of $f$ and $g$.}\label{secd}
The crossing parity $\pi(f,g)$ does not depend on $f$ and $g$
as long as $\phi,\psi:[0;1]\to [0;1]^2$ are continuous and (\ref{f1}) is true.
For showing $\pi(p,q)=1$ we construct special tracks $p$ and $q$.
For tracks $p$ define
$\overline \mathcal L(p):= \bigcup\{l(x,y)\mid x,y\in\mathcal V(p),\ \ x\neq y\}$ ( (\ref{f1b}) and Figure~\ref{fig2}).

\begin{lemma}\label{l14}
There are $5$-approximations $ p = ((s_0,x_0)\ddd(s_k,x_k))$ and
 $q = ((t_0,y_0)\ddd(t_l,y_l))$
of $f$ and $g$, respectively, such that
\begin{eqnarray}
\label{f29}&{\rm card}\ \mathcal V(p)=k+1 \an {\rm card}\ V(q)=l+1 \an\\
\label{f6}&\mathcal V(p)\cap \overline \mathcal L(q)=\emptyset \an
\mathcal V(q)\cap \overline \mathcal L(p)=\emptyset\,.
\end{eqnarray}
\end{lemma}

\proof
There is  a modulus of uniform continuity ${\rm md}$ of $f$ and $g$
(cf. (\ref{f16})).
 There are  some $k$ and rational numbers $s_0:=-1,\ s_1\ddd s_{k-1},s_k:=2$ such that $0<s_{i+1}-s_i<2^{-{\rm md}(5)}$ for
$0\leq i<k$. There are points $x_0\ddd x_k\in\IQ^2$ such that
$\|f(s_0)-x_0\|<2^{-5}$ and $x_i\neq x_j$ for all $i\neq j$.
 Then $ p = ((s_0,x_0)\ddd(s_k,x_k))$  is a $5$-approximation of $f$
 such that ${\rm card}\ \mathcal V(p)=k+1$.
 There are  some $l$ and rational numbers
$t_0:=-1,\ t_1\ddd t_{k-1},t_k:=2$ such that
 $0<t_{i+1}-t_i<2^{-{\rm md}(5)}$ for
$0\leq i<l$.
There is some $y_0\in\IQ^2$ such that $\|g(t_0)-y_0\|<2^{-5}\|$ and
$y_0\not\in\overline\mathcal L(p)$. Suppose $y_0\ddd y_i$  ($i<l$)
have been found.
Then there is some $y_{i+1}\in\IQ^2$ such that
$\|g(s_{i+1})-y_{i+1}\|<2^{-5}$, $y_{i+1}\not\in\{y_0\ddd y_i\}$,
$y_{i+1}\not\in \overline\mathcal L(p)$ and
$(l(y_0,y_{i+1})\cup\ldots \cup l(y_i,y_{i+1}))\cap \mathcal V(p)=\emptyset$.
By induction, we can define $y_0\ddd y_l$. The obviously
(\ref{f29}) and (\ref{f6}) are true for $p$ and $q$.
\qq

\begin{lemma}\label{l11}
$\pi(f,g)=1\,.$
\end{lemma}

\begin{figure}[htbp]
\setlength{\unitlength}{.14ex}
%\linethickness{0.7pt}
\begin{picture}(200,280)(-327,-88)

\put(-180,0){\vector (1,0){460}}
\put(0,100){\line (1,0){100}}
\put(0,-95){\vector (0,1){290}}
\put(100,0){\line (0,1){100}}
\put(0,0){\circle*{6}}
\put(299,0){\makebox(0,0)[cc]{$z_1$}}
\put(-13,190){\makebox(0,0)[cc]{$z_2$}}

\thicklines
\qbezier(0,0)(20,80)(100,100)
\qbezier(0,100)(50,10)(100,0)

\put(-100,0){\line(1,0){100}}
\put(100,0){\line(1,0){100}}
\put(-100,100){\line(1,0){100}}
\put(100,100){\line(1,0){100}}

\put(54,70){\makebox(0,0)[cc]{$\phi$}}
\put(39,29){\makebox(0,0)[cc]{$\psi$}}

\put(-50,90){\makebox(0,0)[cc]{$g$}}
\put(-29,-12){\makebox(0,0)[cc]{$f$}}
\put(125,108){\makebox(0,0)[cc]{$f$}}
\put(150,8){\makebox(0,0)[cc]{$g$}}

\put(8,-12){\makebox(0,0)[cc]{$\phi(0)$}}
\put(8,108){\makebox(0,0)[cc]{$\psi(0)$}}
\put(100,-12){\makebox(0,0)[cc]{$\psi(1)$}}
\put(100,110){\makebox(0,0)[cc]{$\phi(1)$}}

\put(-104,13){\makebox(0,0)[cc]{$f(-1)$}}
\put(-100,110){\makebox(0,0)[cc]{$g(-1)$}}
\put(200,-10){\makebox(0,0)[cc]{$g(2)$}}
\put(200,110){\makebox(0,0)[cc]{$f(2)$}}

\put(-100,0){\circle*{5}}
\put(-100,100){\circle*{5}}
\put(200,0){\circle*{5}}
\put(200,100){\circle*{5}}

\thinlines

\newsavebox{\balle}
\savebox{\balle}{
\qbezier(131,0)(131,54.15)(92.66,92.66)
\qbezier(92.66,92.66)(54.15,131)(0,131)
\qbezier(131,0)(131,-54.15)(92.66,-92.66)
\qbezier(92.66,-92.66)(54.15,-131)(0,-131)
\qbezier(-131,0)(-131,54.15)(-92.66,92.66)
\qbezier(-92.66,92.66)(-54.15,131)(0,131)
\qbezier(-131,0)(-131,-54.15)(-92.66,-92.66)
\qbezier(-92.66,-92.66)(-54.15,-131)(0,-131)
\put(0,0){\circle*{5}}
}
\put(50,50){\usebox{\balle}}

\put(60,45){\makebox(0,0)[cc]{$z$}}

\put(50,0){\line (0,1){50}}
\put(-50,-5){\line(0,1){10}}
\qbezier(-50,0)(50,50)(50,50)
\put(-115,-24){\vector(3,1){61}}
\put(-130,-33){\makebox(0,0)[cc]{$(-0.5,0)$}}

\newsavebox{\rectangle}
\savebox{\rectangle}{
\put(-6.5,-6.5){\line(1,0){69.5}}
\put(-6.5,6.5){\line(1,0){69.5}}
\put(-6.5,-6.5){\line(0,1){13}}
\put(63,-6.5){\line(0,1){13}}
}
\put(-100,0){\usebox{\rectangle}}
\put(146.875,100){\usebox{\rectangle}}

\qbezier(-43,-4)(-43,-4)(147,97)
\put(-43,-4){\circle*{4}}
\put(147,97){\circle*{4}}
\put(-59,-52){\vector(1,3){15}}
\put(163,147){\vector(-1,-3){15}}

\put(-62,-58){\makebox(0,0)[cc]{$x_\beta$}}
\put(164,155){\makebox(0,0)[cc]{$x_\gamma$}}

 \end{picture}
\caption{Intersecting curves $\phi$ and $\psi$ with extensions $f$ and $g$} \label{fig7}
\end{figure}
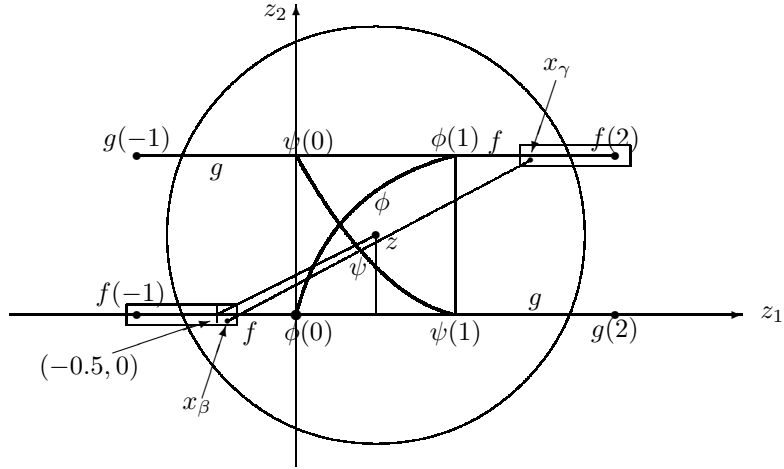

\proof
For $I_0:=J_0:=[-1;2]$ we have
 $f|_{I_0}=f$, $g|_{J_0}=g$ and $\alpha_{I_0J_0}=1$ ((\ref{f8}) and Figure~\ref{fig1}).
Since $2^{-5}< \alpha_{IJ}/16$ by Definition~\ref{d2}, $\pi(f,g)=\pi(p,q)$ for (arbitrary) $5$-approximations $p$ and $q$ of $f$ and $g$ such that (\ref{f29}) and (\ref{f6}). Notice that (\ref{f6}) implies $p\bowtie q$.
In the following we prove $\pi(p,q)=1$.
We define
\begin{eqnarray}
z & := & (0.5,0.5)\,,\\
R_{f0} & := & [-1-2^{-4}; -0.5 +3\cdot 2^{-4}] \times [-2^{-4}; 2^{-4}] \an \\
R_{f1} & := & [1.5-3\cdot 2^{-4}; 2+ 2^{-4}] \times [-2^{-4}; 2^{-4}]\,.
\end{eqnarray}
Figure~\ref{fig7} shows $f$ and $g$, the ball $B(z,1.31)$ and the two boxes $R_{f0}$ and $R_{f1}$.
By the definitions of $f$ and $g$ in (\ref{f15}) and (\ref{f14}),
$f(s)=(s,0)$ for $-1\leq s\leq 0$.
First we prove
\begin{eqnarray}
\label{f28}
h_p(s)\in R_{f0} & \mbox{for} & s<-0.5+2^{-5} \an\\
\label{f17}
\|h_p(s)-z\|< 1.31 & \mbox{for} & -0.5-2^{-5} < s\leq 1.5 + 2^{-5}\an\\
\label{f27}h_p(s)\in R_{f1} & \mbox{for} &
 s > 1.5 -2^{-5}\,.
\end{eqnarray}

We show (\ref{f28}).
Let $s<-0.5+2^{-5}$. Then for some $i$,  $s_i\leq s<s_{i+1}< -0.5 +2^{-4}$.
Then $x_i\in R_{f0}$ since $\|f(s_i)-x_i\|-2^{-5}$. Correspondingly, $x_{i+1}\in R_{f0}$. Since $h_p(s)\in \overline {x_ix_{i+1}}$ and $R_{f0}$ is convex, $h_p(s)=\in R_{f0}$.

(\ref{f27}) is true by symmetry.

\smallskip
We show (\ref{f17}).
If  $ -0.5-2^{-5} <  s\leq -0.5$ then
$\|h_p(s)-z\|\leq \|h_p(s)-f(s)\|+\|f(s)-f(-0.5)\| +\|f(-0.5)-z\|
 <  5\cdot2^{-5} +2^{-5} +\|f(-0.5)-z\|
  \leq 6\cdot 2^{-5}+ \sqrt{ 1+ 1/4})
< 0.19 +  1.12 =1.31$.
If  $ -0.5\leq s\leq 0$ then
$\|h_p(s)-z\|\leq \|h_p(s)-f(s)\|+\|f(s)-z\|<5\cdot 2^{-5} +
\|f(-0.5)-z\| <1.31$ (from above).
Hence $\|h_p(s)-z\|<1.31$  for $ -0.5-2^{-5} < s\leq 0$.
By symmetry, $\|h_p(s)-z\|<1.31$ also for $1\leq s\ <  1.5 + 2^{-5}$.

For $0\leq s\leq 1$, $\|h_p(s)-z\| \leq \|h_p(s)-f(s)\| +\|f(s)-z\| \leq 5\cdot2^{-5}+\sqrt{2}/2 <1.31$.
Therefore  (\ref{f17}) is true.
Since $|s_{i+1}-s_i| <2^{-{\rm md}(5)}\leq 2^{-5}$ there are indices $\beta<\gamma$ such that
\begin{eqnarray}
-0.5-2^{-5}<s_\beta <-0.5+2^{-5} & \an & 1.5-2^{-5}<s_\gamma <1.5+2^{-5}\,.
\end{eqnarray}
 By (\ref{f28}) - (\ref{f27}),
\begin{eqnarray}
\label{f36} x_i \in R_{f0} & \mbox{if} &   i\leq \beta\,,\\
\label{f37}x_i\in  B((z,1.31) & \mbox{if} &  \beta\leq i\leq \gamma\,,\\
\label{f38}x_i \in R_{f1} & \mbox{if} &  \gamma\leq i\,.
\end{eqnarray}
\noindent We simplify the track $p$ step by step by means of Lemma~\ref{l8}.
For $\beta<i\leq\gamma$ define
$$ p_i:=((s_0,x_0)\ddd(s_\beta,x_\beta),(s_i,x_i)\ddd (s_k,x_k))$$
Since the $x_i$ are pairwise different (by ${\rm card}(\mathcal V(p))=k+1$),  $x_\beta\neq x_i$. Therefore, $p_i$ is a track, see Definition~\ref{d1}.
We prove by induction
\begin{eqnarray}\label{f26}
  \pi(p_i,q)=\pi(p,q) & \mbox{for} & \beta<i\leq\gamma\,.
 \end{eqnarray}
Since $p_{\beta+1}=p$,  $ \pi(p_{\beta+1},q)=\pi(p,q) $.
\\
Suppose $ \pi(p_i,q)=\pi(p,q) $ for some $\beta<i<\gamma$. We define
\begin{eqnarray}
p_a & := & ((s_0,x_0)\ddd(s_\beta,x_\beta))\,,\\
p_b & := &  ((s_\beta,x_\beta),(s_i,x_i),(s_{i+1},x_{i+1}))\,,\\
p_c & := &  ((s_\beta,x_\beta),(s_i,z_1),(s_{i+1},x_{i+1}))\\
\nonumber  && \ \ \ \mbox {for some}\ \ z_1\in \overline{x_\beta x_{i+1}}\setminus \{x_\beta,x_{i+1}\}\,,\\
 p_d & := &  ((s_\beta,x_\beta),(s_{i+1},x_{i+1}))\\
 p_e & := & ((s_{i+1},x_{i+1})\ddd(s_k,x_k))
 \end{eqnarray}
We apply Lemma~\ref{l8} to $p_b$ and $p_c$.  We have (cf. (\ref{f30} - \ref{f32}))\\
-- $\{x_\beta,x_{i+1}\}\cap \overline\mathcal L(q) = \emptyset$,
(since $\mathcal V (p)\cap \overline\mathcal L(q)=\emptyset$), \\
-- $\{x_\beta,x_i,x_{i+1},z_2\}\in B(z,1.31)$ (by (\ref{f37})) and\\
-- $\{y_0,y_l\}\cap B(z,1.31)=\emptyset$.\\
(For the last line:  $1.58<\sqrt{1.5^2+0.5^2} =\|g(-1)-z\| \leq \|g(-1)-y_0\| +\|y_0-z\|  <  \|y_0-z\| +2^{-5} <
 \|y_0-z\| + 0.04 $, hence $ \|y_0-z\|> 1.31$. Accordingly
$ \|y_l-z\|> 1.31$.) By Lemma~\ref{l8}, $\pi(p_b,q)=\pi(p_c,q)$.
Since $z_1\in \overline{x_\beta x_{i+1}}\setminus \{x_\beta,x_{i+1}\}$, $\pi(p_c,q)=\pi(p_d,q)$, hence
 $\pi(p_b,q)=\pi(p_d,q)$. Therefore, replacing $p_b$ by $p_d$ in $p_{\beta,i}$ does not change the crossing parity. We show this in detail.

 Since $\mathcal V(p)\cap \overline\mathcal L(q)=\emptyset$, $x_i\not\in \range(h_q)$ for all $i$, hence
\begin{eqnarray*}
{\rm CN}(p_i,q) &=& {\rm CN}(p_a,q) + {\rm CN} (p_b,q) + {\rm CN}(p_e,q) \an\\
{\rm CN}(p_{i+1},q) &= &  {\rm CN}(p_a,q) + {\rm CN} (p_d,q) + {\rm CN}(p_e,q)\,.
\end{eqnarray*}
Since $(m+n)\ \mod\ 2= (m \ \mod \ 2  + n \ \mod \ 2)\ \mod\ 2$,
\begin{eqnarray*}
\pi(p_i,q)\   &=&( \pi(p_a,q) + \pi (p_b,q) + \pi(p_e,q))\ \ \mod \ 2 \\
&=&( \pi(p_a,q) + \pi (p_d,q) + \pi(p_e,q))\ \ \mod \ 2 \\
&=& \pi(p_{i+1},q)\,.
\end{eqnarray*}
By induction
%$\pi(p,q)=\pi(p_\gamma,q)$ where
%
\begin{eqnarray}\label{f39}
\pi(p,q)=\pi(p_\gamma,q)\ \  \mbox{for}&
p_\gamma=((s_0,x_0)\ddd (s_\beta,x_\beta),(s_\gamma,x_\gamma)\ddd(s_k,x_k))\,.
\end{eqnarray}
By (\ref{f28}) and (\ref{f27}), $h_{p_\gamma}$ runs in the box $R_{f1}$ for $-1\leq s\leq s_\beta$, then its graph is the straight line segment $\overline{x_\beta x_\gamma}$  and finally $h_{p_\gamma}$ remains in the box $R_{f_2}$ for $s_\gamma <s$, see Figure~\ref{fig7}.

\begin{figure}[htbp]
\setlength{\unitlength}{.14ex}
%\linethickness{0.7pt}
\begin{picture}(200,210)(-330,-50)

\put(-180,0){\vector (1,0){460}}
\put(0,100){\line (1,0){100}}
\put(0,-60){\vector (0,1){200}}
\put(100,0){\line (0,1){100}}
\put(0,0){\circle*{6}}
\put(299,0){\makebox(0,0)[cc]{$z_1$}}
\put(-13,140){\makebox(0,0)[cc]{$z_2$}}

\put(-100,0){\circle*{5}}
\put(-100,100){\circle*{5}}
\put(200,0){\circle*{5}}
\put(200,100){\circle*{5}}

\newsavebox{\rectanglea}
\savebox{\rectanglea}{
\put(-6.5,-6.5){\line(1,0){69.5}}
\put(-6.5,6.5){\line(1,0){69.5}}
\put(-6.5,-6.5){\line(0,1){13}}
\put(63,-6.5){\line(0,1){13}}
}
\put(-100,0){\usebox{\rectanglea}}
\put(146.875,100){\usebox{\rectanglea}}

\qbezier(-43,-4)(-43,-4)(147,97)
\put(-43,-4){\circle*{4}}
\put(147,97){\circle*{4}}
\put(-59,-52){\vector(1,3){15}}
\put(163,147){\vector(-1,-3){15}}

\put(-62,-58){\makebox(0,0)[cc]{$x_\beta$}}
\put(164,155){\makebox(0,0)[cc]{$x_\gamma$}}

\put(-100,100){\usebox{\rectanglea}}
\put(146.875,0){\usebox{\rectanglea}}

\qbezier(-43,104)(-43,104)(147,3)
\put(-43,104){\circle*{4}}
\put(147,3){\circle*{4}}
\put(-59,152){\vector(1,-3){15}}
\put(163,-47){\vector(-1,3){15}}
\put(-62,158){\makebox(0,0)[cc]{$y_\mu$}}
\put(164,-55){\makebox(0,0)[cc]{$y_\nu$}}

\put(-106,18){\makebox(0,0)[cc]{$R_{f1}$}}
\put(-106,82){\makebox(0,0)[cc]{$R_{g1}$}}
\put(206,18){\makebox(0,0)[cc]{$R_{g2}$}}
\put(206,82){\makebox(0,0)[cc]{$R_{f2}$}}

% urspruenglicher Kreis
%\newsavebox{\ballz}
%\savebox{\ballz}{
%\qbezier(15,0)(15,6.21)(10.61,10.61)
%\qbezier(10.61,10.61)(6.21,15)(0,15)
%\qbezier(15,0)(15,-6.21)(10.61,-10.61)
%\qbezier(10.61,-10.61)(6.21,-15)(0,-15)
%\qbezier(-15,0)(-15,6.21)(-10.61,10.61)
%\qbezier(-10.61,10.61)(-6.21,15)(0,15)
%\qbezier(-15,0)(-15,-6.21)(-10.61,-10.61)
%\qbezier(-10.61,-10.61)(-6.21,-15)(0,-15)
%\put(0,0){\circle*{.5}}
%}
%\put(50,50){\usebox{\ballz}}

 \end{picture}
\caption{The tracks $p_\gamma$ and $q_\nu$.} \label{fig8}
\end{figure}
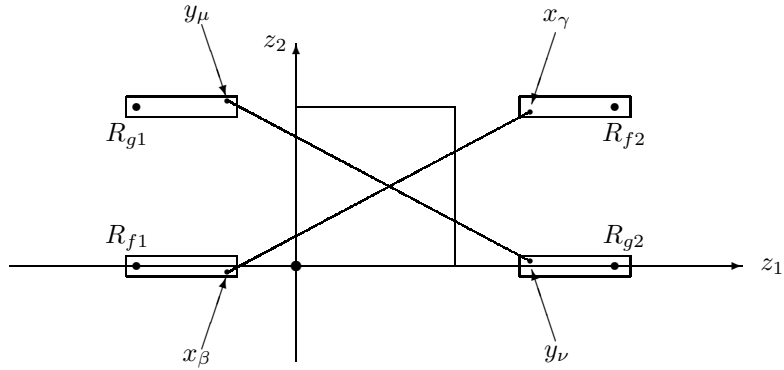

Since $\mathcal V(p_\gamma)\In \mathcal V(p)$,
$\mathcal V(q)\cap \overline \mathcal L(p_\gamma)=\emptyset$.
Keeping $p_\gamma$ fixed we can simplify the track $q$ accordingly.
Therefore there are indices $\mu<\nu$ such that for
$q_\nu:= ((t_0,y_0)\ddd(t_\mu,y_\mu)(t_\nu,y_\nu)\ddd(t_l,y_l))$, $\pi(p_\gamma,q)=\pi(p_\gamma,q_\nu)$, where $h_{q_\nu}$ runs in the box $R_{g1}$ for
$-1\leq t_\mu$, then its graph is the straight line segment $\overline{y_\mu y_\nu}$  and finally $h_{p_\gamma}$ remains in the box $R_{g_2}$ for $t_\nu <t$ as shown in Figure~\ref{fig8}.

Since
$(R_{f1}\cup R_{f2})\cap \range(h_{q_\nu})=\emptyset$ and
$(R_{g1}\cup R_{g2})\cap \range(h_{p_\gamma})=\emptyset$ (we omit the straightforward formal verifications),
 ${\rm CN}(p_\gamma,q_\nu)=1$.
We obtain  $\pi(p,q)=\pi(p_\gamma,q)=\pi(p_\gamma,q_\nu)=1$.
\hfill $\Box$ (Lemma~\ref{l11})

\section{The proof of the main theorem}\label{sece}

For $A\in \IR^2$ and a number $\delta>0$ let
$U(A,\delta):=\{ x\in\IR^2\mid (\exists \,y\in A)\,\|x-y\|<\delta\}$ be the $\delta$-neighborhood of $A$.
From $f$ and $g$ we will compute  sequences $I_0\supseteq I_1\supseteq I_2\supseteq \ldots$ and
$J_0\supseteq J_1\supseteq J_2\supseteq \ldots$ of rational closed intervals such that
$f(\bigcap_i I_i)=g(\bigcap_i J_i)$. First we prove the following lemma.

\begin{lemma}\label{l9}

From $f,g$, rational intervals $I,J\In [-1;2]$ and $n\In \IN$ such that
such that $2^{-n} < \alpha_{IJ}$ and $\pi(f|_I,g|_J)=1$
we can compute a rational interval $K\In I$  such that
\begin{eqnarray}\label{f73}
\alpha_{KJ}>0\,,\
f(K)\In U(g(J),\ee) \an  \pi(f|_K,g|_J)=1\,.
\end{eqnarray}
\end{lemma}

\proof Perform the following computations:

-- From $J$ and $g$ compute the compact set $g(J)$.

-- Compute some sequence $(a_I=s_0,s_1\ddd s_{k-1},s_k=b_I)$ of rational numbers such that $0<s_{i+1}-s_i<2^{-{md}(n+4)}$.\\
 \medskip
For every index $i$ do the following:

-- From $f$ and $s_i$ compute $x_i:=f(s_i)\in \IR^2$.

-- From $x_i$ and $g(J)$ compute $c_i:=d_s(\{x_i\},g(J))\in\IR$.

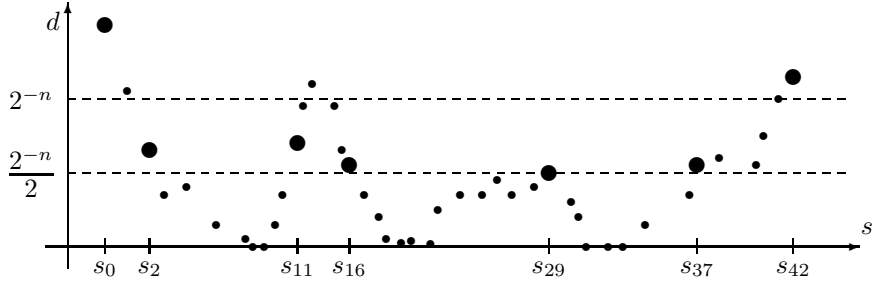
\begin{figure}[htbp]
\setlength{\unitlength}{.65ex}
%\linethickness{0.7pt}
\begin{picture}(120,38)(-30,-3)

\put(-3,0){\vector(1,0){110}}
\put(0,-3){\vector(0,1){36}}
\multiput(0,10)(2,0){53}{\line(1,0){1}}
\multiput(0,20)(2,0){53}{\line(1,0){1}}

\put(-5,10){\makebox(0,0)[cc]{\Large$\frac{\ee}{2}$}}
\put(-5,20){\makebox(0,0)[cc]{$\ee$}}

\put(-2,29){\makebox(0,3)[cc]{$d$}}
\put(108,1){\makebox(0,3)[cc]{$s$}}

\put(5,-1){\line(0,1){2}}
\put(5,-3){\makebox(0,0)[cc]{$s_0$}}
\put(5,30){\circle*{2}}
\put(8,21){\circle*{1}}
\put(11,-1){\line(0,1){2}}
\put(11,-3){\makebox(0,0)[cc]{$s_2$}}
\put(11,13){\circle*{2}}
\put(13,7){\circle*{1}}
\put(16,8){\circle*{1}}
\put(20,3){\circle*{1}}
\put(24,1){\circle*{1}}
\put(25,0){\circle*{1}}
\put(26.5,0){\circle*{1}}
\put(28,3){\circle*{1}}
\put(29,7){\circle*{1}}
\put(31,-1){\line(0,1){2}}
\put(31,-3){\makebox(0,0)[cc]{$s_{11}$}}
\put(31,14){\circle*{2}}
\put(31.8,19){\circle*{1}}
\put(33,22){\circle*{1}}
\put(36,19){\circle*{1}}
\put(37,13){\circle*{1}}
\put(38,-1){\line(0,1){2}}
\put(38,-3){\makebox(0,0)[cc]{$s_{16}$}}
\put(38,11){\circle*{2}}
\put(40,7){\circle*{1}}
\put(42,4){\circle*{1}}
\put(43,1){\circle*{1}}
\put(45,.5){\circle*{1}}
\put(46.4,.7){\circle*{1}}
\put(49,.4){\circle*{1}}
\put(50,5){\circle*{1}}
\put(53,7){\circle*{1}}
\put(56,7){\circle*{1}}

%\put(58,-1){\line(0,1){2}}
%\put(58,-3){\makebox(0,0)[cc]{$s_{26}$}}
\put(58,9){\circle*{1}}
\put(60,7){\circle*{1}}
\put(63,8){\circle*{1}}

\put(65,-1){\line(0,1){2}}
\put(65,-3){\makebox(0,0)[cc]{$s_{29}$}}
\put(65,10){\circle*{2}}
\put(68,6){\circle*{1}}
\put(69,4){\circle*{1}}
\put(70,0){\circle*{1}}
\put(73,0){\circle*{1}}
\put(75,0){\circle*{1}}
\put(78,3){\circle*{1}}
\put(84,7){\circle*{1}}
\put(85,-1){\line(0,1){2}}
\put(85,-3){\makebox(0,0)[cc]{$s_{37}$}}
\put(85,11){\circle*{2}}
\put(88,12){\circle*{1}}
\put(93,11){\circle*{1}}
\put(94,15){\circle*{1}}
\put(96,20){\circle*{1}}
\put(98,-1){\line(0,1){2}}
\put(98,-3){\makebox(0,0)[cc]{$s_{42}$}}
\put(98,23){\circle*{2}}

\end{picture}
\caption{$d_i\approx d_s( f(s_i),g(J))$} \label{fig11}
\end{figure}

-- From $c_i\in\IR$ and $n$ compute some $d_i\in\IQ$ such that
 $|c_i-d_i|< \ee /16$.\\
Figure~\ref{fig11} shows the points $d_i\approx d_s(f(s_i),g(J))$.

Suppose $|s_i-s|<2^{-{md}(n+4)}$   and $d_i\geq \ee/2$.
Then for every $z\in g(J)$, \\
$ \ee/2\leq  d_i < c_i +  \ee/16\leq \|x_i-z\|+ \ee/16
\leq \|f(s_i)-f(s)\| + \|f(s)-z\|+ \ee/16 <
 \ee/16 + \|f(s)-z\|+ \ee/16 $, hence $\|f(s)-z\|>6/16\eee$.
Therefore,
\begin{eqnarray}
\nonumber &\mbox{if  $|s_i-s|<2^{-{md}(n+4)}$ and $d_i\geq \ee/2$ then}\\
&\label{f46}
(\forall\;z\in f(J))\|f(s)-z\|>6/16\eee \,.
\end{eqnarray}

Suppose $|s_i-s|<2^{-{md}(n+4)}$ and $d_i < \ee/2$.
Then $c_i<d_i+\ee/16\leq 9/16\eee$, hence for some $z\in g(J)$,
$\|f(s_i)-z\| < 9/16\eee$ and
$\|f(s)-z\|\leq \|f(s)-f(s_i)\| +\|f(s_i)-z\|\leq 10/16 \eee$.
Therefore,
\begin{eqnarray}
\nonumber&\mbox{ if $|s_i-s|<2^{-{md}(n+4)}$ and $d_i < \ee/2$ then}\\
\label{f48}
&(\exists\;z\in f(J)) \|f(s)-z\|<10/16\eee \,.
\end{eqnarray}
Since $2^{-n} <\alpha_{IJ}$ and $s_0=a_I$, $\ee <  d_s(f(a_I),g(J))=c_0<d_0+\ee/16$, hence
\begin{eqnarray}\label{f87}
d_0>\ee/2, & \mbox{correspondingly,} &d_k>\ee/2\,.
\end{eqnarray}
Therefore, we can find a sub-sequence $0=i_0<i_1<\ldots <i_m=k$
of $0,1 \ddd k$ such that
\begin{eqnarray}
\label{f70}
&(\forall\, j)\,d_{i_j}\geq \ee/2 \an\\
\label{f55}
 \mbox{either} & (\forall \, i_j<l<i_{j+1})\, d_l\geq  \ee/2 & \\
\label{f69} \mbox{or} &
(\forall \, i_j<l<i_{j+1})\, d_l <  \ee/2\,.
\end{eqnarray}
Let $K_j:=[s_{i_j};s_{i_{j+1}}]$. (In Figure~\ref{fig11}, $m=6$ and $(i_0,i_1\ddd i_6)=(0,2,,11,16,29,37,42)$.)

Suppose (\ref{f55}), that is, $(\forall \, i_j<l<i_{j+1})\, d_l\geq  \ee/2$. (In Figure~\ref{f11} this is true for $j=0,2,5$.)
 By (\ref{f46}) for all $z\in g(J)$ and $s\in K_j$, $\|f(s)-z\|>6/16\eee$,
hence $f(K_j)\cap g(J)=\emptyset$, therefore,
\begin{eqnarray}\label{f71}
\pi(f|_{K_j},g|_{J})=0\,.
\end{eqnarray}
Suppose (\ref{f69}), that is, $(\forall \, i_j<l<i_{j+1})\, d_l <  \ee/2$. (In Figure~\ref{f11} this is true for $j=1,3,4$.)
By (\ref{f48})
for all $s\in K_j$ there is some $z\in g(J)$ such that
$\| f(s)-z\|<10/16\eee$, hence $f(K_j)\In U(g(J), \ee)$.
Furthermore by(\ref{f70}), $\ee/2-\ee/16\leq d_{i_j}-\ee/16<c_{i_j}=d_s(\{f(s_{i_j})\},g(J))$. The same is true for
$(j+1)$ instead of $j$, therefore
$0<d_s(\{f(s_{i_j}),  f(s_{i_{j+1}})  \},g(J))$.

Since $K_j\In I$ and $\alpha_{IJ}>0$,
$d_s(f(K_j),\{g(a_J),g(b_J)\}\geq d_s(f(I),\{g(a_J),g(b_J)\}\geq \alpha_{IJ}$. In summary,
\begin{eqnarray}\label{f72}
  0<\alpha_{K_jJ} \an  f(K_j)\In U(g(J), \ee)\,.
\end{eqnarray}
By Lemma~\ref{l5},
$$ \big(\pi(f|_{K_0}, g|_{J})
+\cdots + (\pi(f|_{K_{m-1}}, g|_{J})\big) \ \mod\;2
=\pi(f|_I, g|_{J})= 1\,.$$
Therefore, there is some $j$ such that $\pi(f|_{K_j}, g|_{J})=1$.
By (\ref{f71}), (\ref{f55}) cannot be true for $j$, hence (\ref{f69}) is true for $j$, hence (\ref{f72}) is true for $K_j$.

For computing $K$ from $f,g,I,J$ and $n$, first compute $(d_0 \ddd
d_k)$, then compute the sub-sequence $0=i_0<i_1<\ddd<i_m=k$
of $0,1 \ddd k$ such that (\ref{f70}--\ref{f69}). The find some $j$ such that $\pi(f|_{K_j}, g|_{J})=1$. Let $K:=K_j$. Then (\ref{f73})  is true by (\ref{f72}).
\qq

\begin{lemma}\label{l10}
From functions $f$ and $g$, closed rational intervals $I,J\In [-1;2]$ and $m\in\IN$ such that $\alpha_{IJ}>0$ and $\pi(f|_I,g|_J)=1$ we can compute closed rational intervals $I',J'$ such that
\begin{eqnarray}
\label{f88}& \alpha_{I'J'}>0  \an \pi(f|_{I'},g|_{J'})=1\,,\\
\label{f89}&I'\In I \an  J'\In J\,, \\
\label{f90}&f(I')\In  U(g(J),2^{-m})\an g(J')\In U(f(I'),2^{-m})\,.
\end{eqnarray}
\end{lemma}

\proof
Compute $\alpha_{IJ}$ and some number $n> m$ such that
$2^{-n} < \alpha_{IJ}$.

By Lemma~\ref{l9} we can compute a rational interval $I'\In I$ ($I_m=K$ in the lemma) such that by (\ref{f73}),
\begin{eqnarray}
\label{f74}\alpha_{I'J}>0,\ \ f(I')\in U(g(J),2^{-n}) \an \pi(f|_{I'},g|_{J})=1\,.
\end{eqnarray}
Again by Lemma~\ref{l9} with ($f$ and $g$ exchanged) first compute
$ \alpha_{I'J}$ and
some number $n'>m $ such that $2^{-n'}<  \alpha_{I'J}$, then compute some $J'\In J$ such that
\begin{eqnarray}
\label{f75}\alpha_{I'J'}>0,\ \ g(J')\in U(f(I'),2^{-n'}) \an
\pi(f|_{I'},g|_{J'})=1\,.
\end{eqnarray}
(\ref{f88}) follows from (\ref{f75}), (\ref{f89}) is true by the construction and (\ref{f90}) follows from  (\ref{f74}) and (\ref{f75}) by $m<n$ and $m<n'$.
\qq

\begin{lemma}\label{l12}
From functions $f$ and $g$ (as defined in (\ref{f15}) and (\ref{f14})) we can compute sequences $I_0\supseteq I_1,\supseteq I_2,\ldots$ and $J_0\supseteq J_1\supseteq J_2,\ldots$ of rational closed intervals such that $f(\bigcap_m I_m)=g(\bigcap_m J_m)$.
\end{lemma}

\proof By Lemma~\ref{l10} starting with $I_0:=J_0:=[-1;2]$ we can compute sequences $I_0\supseteq I_1,\supseteq I_2,\ldots$ and $J_0\supseteq J_1\supseteq J_2,\ldots$ of rational closed intervals such that for all $m\geq 1$,
\begin{eqnarray}
\label{f76}
f(I_m)\In  U (g(J_{m-1}),2^{-(m-1)}) \an g(J_m)\In
 U(f(I_m)),2^{-m})\,.
\end{eqnarray}
The following equations follow from continuity of $f$ and $g$.
\begin{eqnarray}
\label{f77}
&& f(\bigcap_m I_m)=\bigcap_m f(I_m) =\bigcap_m U(f(I_m),2^{-m})\,,\\
\label{f78}
&& g(\bigcap_m J_m)=\bigcap_m g(J_m) =\bigcap_m U(g(J_m),2^{-m})\,.
\end{eqnarray}
We give elementary proofs. Since $f(\bigcap_kI_k) \In f(I_m)\In U(f(I_m),2^{-m})$ for all $m$,
$f(\bigcap_k I_k)\In \bigcap_m f(I_m) \In \bigcap_m U(f(I_m),2^{-m})$.

Suppose $x\in \bigcap_m f(I_m)$. Then there is a sequence $a_0,a_1,a_2,\ldots$ such that for all $m$, $a_m\in I_m$ and $f(a_m)=x$.
Since $I_0$ is compact there is a subsequence $a_{m_0},a_{m_1},a_{m_2},\ldots$ converging to some $a\in I_0$. Since $f$ is continuous, $x=\lim_if(a_{m_i}) =f(\lim_i a_{m_i})=f(a)$. Assume $a\not\in I_{m_k}$ for some $k$. Then for some $j>k$,
$a_{m_j}\not\in I_{m_k}$, hence  $a_{m_j}\not\in I_{m_j}$, a contradiction.
Therefore $a\in I_{m_k}$  for all $k$, hence
$a\in \bigcap_m I_m$ and $x=f(a)\in f(\bigcap_m I_m)$.

Suppose $x\in \bigcap_m U(f(I_m),2^{-m})$. Assume $x\not\in f(I_k)$ for some $k$. Since $f(I_k)$ is compact $x\not\in U(f(I_k),2^{-j})$ for some $j>k$. Therefore, $x\not\in  U(f(I_j),2^{-j})$, hence
$x\not \in \bigcap_m U(f(I_m),2^{-m})$, a contradiction. Therefore
$x\in f(I_k)$ for all $k$, hence $x\in \bigcap_m f(I_m)$.

Therefore (\ref{f77}) is true. By symmetry also (\ref{f78}) is true.
By (\ref{f76}--\ref{f78}),

 $ f(\bigcap_m I_m)=\bigcap_m f(I_m)=\bigcap_{m\geq 1} f(I_m)
    \In \bigcap_{m\geq 1} U (g(J_{m-1}),2^{-(m-1)})
=\bigcap_m U (g(J_m),2^{-m})=g(\bigcap_mJ_m)$. Accordingly,
$g(\bigcap_mJ_m) \In  f(\bigcap_m I_m) $.
\qq\\

After these preparations Theorem~\ref{t2} can be proved straightforwardly. For the canonical representation (most conveniently the canonical multi-representation \cite{WG09}) of functions $h\pf \IR\to\IR^2$ the function $f$ can be computed from $\phi$ and the function $g$ can be computed from $\psi$. By the outer representation of the set $\mathcal I$  of the closed real intervals, a name of $S$ is a sequence  $I_0\supseteq I_1\supseteq\ldots$ of closed rational intervals such that
$S=\bigcap _n I_n$. Therefore, by Lemma~\ref{l12} the multi-function
$(\phi,\psi) \mmto (S_\phi,S_\psi)$ such that $\phi(S_\phi)=\psi(S_\psi)$ is computable. Thus we have proved Theorem~\ref{t2}.

The main theorem from \cite{Wei19} follows straightforwardly from Theorem~\ref{t2}.

\newpage

\begin{corollary}\label{cor2}$ $
\begin{enumerate}
\item  If $\phi$ and $\psi$ in Theorem~\ref{t1} are computable then there are computable numbers $a$ and $b$ such that $\phi(a)\in \range(\psi)$ (hence $x:=\phi(a)\in \range(\phi)\cap  \range(\psi)$)  and
    $\psi(b)\in\range(\phi)$ (hence $y:=\psi(b)\in \range(\phi)\cap  \range(\psi)$).

\item Restricted to the pairs $(\phi,\psi)$ which have a unique intersection point, the point $x$ such that
    $\{x\}=\range(\phi)\cap  \range(\psi)$
    can be computed from $\phi$ and $\psi$.
\item Restricted to the pairs  $(\phi,\psi)$  such that $\phi(a)=\psi(b)$ for a unique pair $(a,b)\in[0;1]^2$
    the function $(\phi,\psi)\mapsto (a,b)$ is computable.

\end{enumerate}
\end{corollary}

\proof

1. Suppose $S_\phi$ has length $>0$. Then $a\in S_\phi$ for some $a\in\IQ$, which is a computable real number such that
$\phi(a)\in \psi(S_\beta)\In \range(\psi)$.

Suppose $S_\phi$ has length $0$, that is, $S_\phi=\{a\}$ for some $a$. Since  $\{ a\} =S_\phi=\bigcap_n I_n$ for a computable sequence $I_0\supseteq I_1\supseteq\ldots$ of rational intervals, $a$ is a computable real number such that $\phi(a)\in \psi(S_\beta)\In\range(\psi)$.

By symmetry there is a computable number $b$ such that $\psi(b)\in\range(\phi)$.

\medskip
 2. Notice that possibly $S_\phi$ and $S_\psi$ have positive lengths. Suppose $\phi$ and $\psi$ have a single intersection point  $x$. Then $\{x\}=\phi(S_\phi)=\psi(S_\psi)$.
For every $m>0$, $\emptyset =\{x\} \setminus B(x,2^{-m})=\bigcap_{n\in\IN}\phi(I_n)\setminus B(x,2^{-m}) =
\bigcap_{n\in\IN}(\phi(I_n)\setminus B(x,2^{-m}))$ by (\ref{f77}). The countable intersection of closed subsets of the compact set $[0;1]^2$ is empty. Therefore finitely many suffice: there is some $N$ such that
$\emptyset = \bigcap_{n\leq N}(\phi(I_n)\setminus B(x,2^{-m}))=
\bigcap_{n\leq N}\phi(I_n)\setminus B(x,2^{-m})=\phi(I_N)\setminus B(x,2^{-m})$,  hence $\phi(I_N)\In B(x,2^{-m})$.
Therefore, for every $m>0$ there are some $N\in\IN$ and some $z\in\IQ^2$ such that $x\in \phi(I_N)\In B(z,2\cdot 2^{-m})$.
Since $(\phi,N)\mapsto \phi(I_N)$ is computable and the set
$(K,z,j)$ such that $K$ is compact and $K\In B(z,j)$ is c.e. \footnote{computably  enumerable or recursively enumerable}  \cite{WG09}, from $\phi$ and the list $I_0,I_1,\ldots$ we can compute a sequence of rational balls contracting to $x$. Therefore we can compute the single intersection point of $\phi$ and $\psi$.
\medskip

 3. We can compute a sequence $I_0\supseteq I_1\supseteq \ldots$ converging to $a$. Therefore we can compute $a$. Correspondingly we can compute $b$.
\qq\\

In Corollary~\ref{cor2}.1 in general $\phi(a)\neq\psi(b)$. From $\phi$ and $\psi$ we cannot compute some $a$ such that $\phi(a)\in\range(\psi)$ or some $x\in \range(\phi)
\cap \range(\psi)$.
We do not even know whether for the computable number $a$ there is a computable number $c$ such that $\phi(a) = \psi(c)$.
In Corollary~\ref{cor2}.2 from $\phi$ and
$\psi$ we cannot compute some $a$  such that $\phi(a)=x$.

\bibliographystyle{plain}

%\bibliography{cca-2019-11-03,meinebib_2019-03-12}

\begin{thebibliography}{10}

\bibitem{BHW08}
Vasco Brattka, Peter Hertling, and Klaus Weihrauch.
\newblock A tutorial on computable analysis.
\newblock In S.~Barry Cooper, Benedikt L\"owe, and Andrea Sorbi, editors, {\em
  New Computational Paradigms: Changing Conceptions of What is Computable},
  pages 425--491. Springer, New York, 2008.

\bibitem{Grz55}
Andrzej Grzegorczyk.
\newblock Computable functionals.
\newblock {\em Fundamenta Mathematicae}, 42:168--202, 1955.

\bibitem{Grz57}
Andrzej Grzegorczyk.
\newblock On the definitions of computable real continuous functions.
\newblock {\em Fundamenta Mathematicae}, 44:61--71, 1957.

\bibitem{Kus84}
Boris~Abramovich Ku{\v{s}}ner.
\newblock {\em Lectures on Constructive Mathematical Analysis}, volume~60 of
  {\em Translations of Mathematical Monographs}.
\newblock American Mathematical Society, Providence, 1984.

\bibitem{Kus99}
Boris~Abramovich Ku{\v{s}}ner.
\newblock Markov's constructive analysis; a participant's view.
\newblock {\em Theoretical Computer Science}, 219:267--285, 1999.

\bibitem{Lac55a}
Daniel Lacombe.
\newblock Extension de la notion de fonction r\'{e}cursive aux fonctions d'une
  ou plusieurs variables r\'{e}elles {I}.
\newblock {\em Comptes Rendus Acad\'emie des Sciences Paris}, 240:2478--2480,
  June 1955.
\newblock Th\'{e}orie des fonctions.

\bibitem{Man76a}
S.~N. Manukyan.
\newblock O nekotorykh topologicheskikh osobennostyakh konstruktivnykh prostykh
  dug. (in {R}ussian).
\newblock In B.A. Kushner and A.A. Markov, editors, {\em Issledovaniya po
  teorii algorifmov i matematicheskoy logike}, volume~2, pages 122--129.
  Vychislitel'ny Tsentr AN SSSR, Moscow, 1976.

\bibitem{TW11}
Nazanin~R.\ Tavana and Klaus Weihrauch.
\newblock Turing machines on represented sets, a model of computation for
  analysis.
\newblock {\em Logical Methods in Computer Science}, 7(2):2:19, 21, 2011.

\bibitem{Wei00}
Klaus Weihrauch.
\newblock {\em Computable Analysis}.
\newblock Springer, Berlin, 2000.

\bibitem{Wei19}
Klaus Weihrauch.
\newblock Computable planar curves intersect in a computable point.
\newblock {\em Computability}, 8(3, 4):399--415, 2019.

\bibitem{WG09}
Klaus Weihrauch and Tanja Grubba.
\newblock Elementary computable topology.
\newblock {\em J.UCS}, 15(6):1381--1422, 2009.

\end{thebibliography}
\end{document}